\documentclass[12pt]{article}
\usepackage{graphicx}
\usepackage{amssymb}
\usepackage{natbib}
\usepackage{url}
\usepackage{multirow}
\usepackage{booktabs}
\usepackage{amsthm, amsmath, amsfonts, amssymb,mathrsfs}
\usepackage{enumerate}
\usepackage{subcaption}
\usepackage{amsthm}

\usepackage[colorlinks,linkcolor=blue,anchorcolor=black,citecolor=blue,urlcolor=black]{hyperref}
\theoremstyle{plain}
\newtheorem{theorem}{Theorem}
\newtheorem{lemma}{Lemma}

\theoremstyle{definition}
\newtheorem{assumption}{Condition}

\newtheorem{example}{Example}

\usepackage[plain,ruled,noend]{algorithm2e}
\usepackage{multirow, multicol}
\usepackage{float, subcaption}
\usepackage{bm, enumerate, xcolor}
\newcommand{\red}{\color{red}  }
\newcommand{\blue}{\color{blue}  }

\graphicspath{{figures/}}

\def\T{{ \mathrm{\scriptscriptstyle T} }}

\newcommand{\E}{\mathbb{E}}

\newcommand{\loss}{\mathscr{L}}
\newcommand{\risk}{\mathcal{R}}

\newcommand{\data}{\mathcal{D}}
\newcommand{\cali}{\mathrm{cal}}
\newcommand{\train}{\mathrm{train}}
\newcommand{\region}{\mathcal{A}}
\newcommand{\sparse}{\mathrm{sparse}}
\newcommand{\dense}{\mathrm{dense}}
\newcommand{\s}{\mathrm{score}}

\newcommand{\lambdaCali}{\widehat{\lambda}^*_{n_\cali}}
\newcommand{\lambdaCRC}{\widehat{\lambda}^\diamond}
\newcommand{\lambdaExcess}{\widehat{\lambda}_{n_{\mathrm{cal}}}^\dagger}
\newcommand{\blind}{1}

\begin{document}

\def\spacingset#1{\renewcommand{\baselinestretch}%
{#1}\small\normalsize} \spacingset{1}

\newenvironment{shrinkeq}[1]%
{ \bgroup
  \addtolength\abovedisplayshortskip{#1}
  \addtolength\abovedisplayskip{#1}
  \addtolength\belowdisplayshortskip{#1}
  \addtolength\belowdisplayskip{#1}}
{\egroup\ignorespacesafterend}


\if1\blind
{
  \title{\bf Trustworthy Decisions in Reliability Set Estimation under Insufficient Model Information}
  \author{
   Holger Dette$^1$, Zhengfu Liu$^2$, and Jun Yu$^2$\thanks{
  		All the authors have equally contributed to this work. Authors are listed alphabetically by their surname. The corresponding author is Jun Yu (yujunbeta@bit.edu.cn).  
        The work of  HD  has been partially supported
by the Deutsche Forschungsgemeinschaft (DFG):
TRR 391 {\it Spatio-temporal Statistics for the Transition of Energy and Transport} (520388526);
Research unit 5381 \textit{Mathematical Statistics in the Information Age} (460867398).}\hspace{.2cm}\\
  	1: Ruhr-Universit\"{a}t Bochum, Fakult\"{a}t f\"{u}r Mathematik\\ 44780 Bochum, Germany\\
    2: School of Mathematics and Statistics, Beijing Institute of Technology\\100081 Beijing, China\\
   }
  \maketitle
} \fi

\if0\blind
{
  \bigskip
  \bigskip
  \bigskip
  \begin{center}
    {\LARGE\bf Trustworthy Decisions in Reliability Set Estimation under Insufficient Model Information }
\end{center}
  \medskip
} \fi

\bigskip
\begin{abstract}

Reliability set estimation identifies input regions where a response probability exceeds a target level, bridging estimation and safety-critical decisions. Practitioners typically start with a working model, an imperfect approximation of the true response surface. Relying on this imperfect model may incur decision risk, potentially certifying unsafe regions as safe. We develop a unified framework that turns such a working model into a trustworthy decision rule. First, a modeling-then-calibration procedure decouples estimation from decision. Since the true set is unobservable, we introduce an asymmetric, observable surrogate loss and use a separate calibration set to select a bias-correcting threshold, reducing decision risk and achieving $O_P(1/n)$ volume convergence. Second, we leverage conformal risk control with the surrogate loss to control false inclusion risk, which is the most safety-critical error, at a pre-specified level regardless of working model quality. Together, these calibration procedures show that a separate calibration set is necessary for risk control. Third, an adaptive design concentrates observations on the reliability set and its boundary, improving model quality where errors most affect decisions while controlling budget elsewhere. 
Numerical studies show not only more accurate set estimates but also calibrated finite-sample risk control that classical plug-in methods lack.

\end{abstract}

\noindent%
\vfill

\newpage
\spacingset{1.75} 
\section{Introduction}\label{sec:intro}
\vspace{-0.25cm}

Reliability set estimation bridges statistical estimation and safety‑critical decision making.
For a binary outcome and a user‑specified probability threshold, the reliability set is the collection of input points where the success probability exceeds that threshold.
Its estimate defines the region where an engineer can safely operate \citep{movahedian2019les,Novik2022,NovikChristensen2024}, which dose range yields a minimal efficiency~\citep{bretz2006,dette2010optimal,frieri2023design}, or which material composition passes a safety standard~\citep{Jensen2020,Michalchuk2021,Cumming2024}.
In each of these settings, an incorrect reliability assessment can have severe consequences.
Although such decision problems originated in the mid‑20th century with the analysis of sensitivity experiments~\citep{dixon1948method}, their relevance has grown substantially in recent years 
\citep{Christensen2024a,rougier2024bayesian}.
The practical stakes are high, yet the central statistical obstacle is fundamental: the true success probability function is never known, but reliable decisions must still be made.

Faced with this obstacle, the practitioner turns to a working model: a useful approximation of the true response surface that encodes domain knowledge, guides data collection, and gives shape to the reliability set. 
Yet any such model is, by definition, a simplification, and its discrepancies from the truth are rarely benign. 
The plug-in approach \citep{langlie1963reliability,young1994estimation} trusts the model literally, converting a fitted surface directly into a decision boundary. 
Any bias in the model, whether from an overly rigid parametric form or from insufficient data to fit a richer one, becomes bias in the decision. 
Thus, the classical plug-in approach fails to translate the estimation into a decision with controlled risk. 
The practitioner is therefore left with no bound on the decision error, especially in finite samples and with insufficient knowledge of model specification, and existing methodology offers no mechanism to provide one.

An alternative strategy, stochastic approximation~\citep{robbins1951stochastic,LAI1985asymptotic,joseph2004efficient}, sidesteps global modeling entirely and targets the decision boundary as a root‑finding problem. 
By avoiding model specification altogether, it escapes the possible bias of the plug‑in approach. 
But it inherits a different limitation: as the input dimension grows, the decision boundary becomes a continuous hypersurface whose exhaustive exploration is prohibitively expensive, and sequential sampling converges slowly to this surface, often missing entire portions of the reliability set. 
More fundamentally, like the plug‑in approach, stochastic approximation offers no finite‑sample mechanism to bound the resulting decision error.

The practitioner is thus faced with two unsatisfactory options: build a model and accept that model misspecification may translate into uncontrolled decision risk, or avoid modeling and accept that the reliability set is difficult to map reliably in high dimensions. What is missing is a framework that can exploit a simple, low-variance working model without requiring it to be correctly specified, while providing an explicit finite-sample guarantee for the resulting decision.

\textbf{Our Contributions.} 
To address this gap, we develop a unified framework that turns an imperfect working model into a trustworthy decision rule.
The framework can be summarized as modeling then calibration: we use calibration to make the final decision trustworthy, and we use adaptive design to improve the quality of the modeling stage, which in turn strengthens the calibration step. 
\begin{itemize}
    \item[(I)] We decouple model estimation from
threshold calibration: the working model captures the overall shape of the response surface, and an independent calibration set selects the threshold that corrects systematic decision errors. 
Theorem~\ref{lem:consistency} proves that calibration never increases risk, relative
to the fixed-threshold plug‑in rule, and that under a mild conditions on the working model the calibrated estimator recovers the true reliability set exactly. 
For a general working model, the volume of the symmetric difference converges at a fast parametric rate, $O_P(n^{-1})$, in the calibration sample size (Theorem~\ref{thm:calibration}).
Calibration is therefore a low-cost safeguard: a modest calibration set suffices even when the model is imperfect, and when the model happens to be correctly specified, calibration incurs no asymptotic efficiency loss (Theorems~S.1 and ~S.2 in the Supplement).
    \item[(II)] Controlling the total risk is not enough: in safety-critical applications, the practitioner needs an explicit bound on the most consequential error---falsely certifying an unsafe input as reliable.
We extend conformal risk control to reliability set estimation. 
The obstacle is that the true set membership is unobservable, so standard losses cannot be evaluated. We overcome this by applying conformal risk control to the proposed surrogate loss in the calibration step, yielding an upper bound on the false inclusion risk that holds regardless of model quality, and a matching lower bound that shows the guarantee is not overly conservative (Theorems~\ref{thm:validity}--\ref{thm:lower_bound}).
Moreover, we show that the proposed practical method controls the false inclusion risk asymptotically to the target level under mild conditions (Theorem~\ref{thm:validity1}). 
The result is a model-agnostic safety certificate: the chance of falsely certifying an unsafe region as reliable is explicitly bounded at a user-specified level, without assuming the working model to be correct (or consistent).
   
    \item[(III)]  The quality of the working model ultimately determines how much calibration must compensate. To improve model quality without wasting experimental resources, we propose a design criterion that concentrates observations where they matter for the final decision. The criterion has two requirements. First, the design should be space-filling on a critical region covering the reliability set and its boundary; this is minimax-optimal for response surface estimation over the reliability set. 
Second, the fraction of the budget spent outside this region should be explicitly controlled. We implement these
principles via an adaptive algorithm that learns the critical region online, without prior knowledge of the target set, and we prove a bound on the expected proportion of misallocated observations (Theorem~\ref{thm:error_count_bound}). 
This improves model quality where it matters for the final decision.

\end{itemize}

A vehicle safety study calibrated to real crash data confirms this: under a tight experimental budget, our approach reduces the symmetric difference by over 65\% compared to the classical plug‑in benchmark while strictly controlling the false inclusion risk at the nominal level.

\section{A New Approach for Reliability Set Estimation}\label{sec:method}

For  $\bm{x} \in [0,1]^d$, let  $Y \in \{0,1\}$ follow a Bernoulli distribution with success probability $\Pr(Y=1 \mid \bm{x})$.  
For a pre‑specified threshold $p_R \in (0,1)$, the reliability set is defined by
\begin{equation}\label{eq:reliable_set}
    \mathcal{X}_R := \big\{ \bm{x} \in [0,1]^d \mid \Pr(Y=1 \mid \bm{x}) \ge p_R \big\}.
\end{equation}
Our goal is to estimate $\mathcal{X}_R$ from data.  
We assume both $\mathcal{X}_R$ and its complement have positive Lebesgue measure, and, for mathematical convenience, that the boundary $\{\bm{x}\mid \Pr(Y=1 \mid \bm{x}) = p_R\}$ has Lebesgue measure zero.

The classical approach starts from a working model $\pi(\bm{x})$ for $\Pr(Y=1 \mid \bm{x})$ and estimates the reliability set via the plug‑in rule
$\hat{\mathcal{X}}_R = \{\bm{x}\mid \pi(\bm{x}) \ge p_R\}$.
In practice, however, $\pi$ is at best a useful but imperfect approximation, which may be derived from an under‑parameterized model, engineering heuristics, or limited training data. 
Under such insufficient model information, estimation errors in $\pi$ near the decision boundary $\{\bm{x}\mid \Pr(Y=1 \mid \bm{x}) = p_R\}$
are directly inherited by the set estimator, silently turning model bias into decision errors.

To make this precise, we quantify the cost of a set estimator
$\hat{\mathcal{X}}_R$ by the following risk
\begin{equation}\label{eq:risk11}
    \risk(\hat{\mathcal{X}}_R) := \E\Bigl\{
        (p_R - Y) \, \mathbb{I}\bigl(\bm{x} \in \hat{\mathcal{X}}_R,\,
        \bm{x} \notin \mathcal{X}_R\bigr)
        + (Y-p_R) \, \mathbb{I}\bigl(\bm{x} \not\in \hat{\mathcal{X}}_R,\,
        \bm{x} \in \mathcal{X}_R\bigr)
    \Bigr\}.
\end{equation}
The first term penalizes false inclusion while the second penalizes false exclusion.
The risk $\risk(\hat{\mathcal{X}}_R)$ is non‑negative and vanishes if and only if
$\hat{\mathcal{X}}_R = \mathcal{X}_R$ almost surely.
When $\pi$ equals the true conditional probability, the classical plug-in estimator is optimal.
If $\pi$ is misspecified, the fixed threshold $p_R$ offers no protection,
and the resulting risk can be arbitrarily far from the achievable minimum.

To reduce the risk incurred by an imperfect working model, we generalize the threshold $p_R$ as a tunable parameter.
Specifically, we consider a family of set estimators
$\hat{\mathcal{X}}_R(\lambda) = \{\bm{x}\mid \pi(\bm{x}) \ge 1-\lambda\}$ indexed by a scalar $\lambda$, where a smaller $\lambda$ yields a more conservative set  and  the classical plug‑in corresponds to the choice $\lambda = 1-p_R$.
With this parameterization, the risk becomes a function of $\lambda$ (therefore, we use the notation $\risk(\lambda):= \risk(\hat{\mathcal{X}}_R(\lambda))$ for short), and we select $\lambda$ to minimize risk \eqref{eq:risk11} for the final decision. 
In this section, we will develop this calibration procedure to reduce the risk relative to the classical approach.
Beyond risk minimization, the calibration set can also be used to provide a finite‑sample, model‑agnostic guarantee on the false inclusion risk (the first part in \eqref{eq:risk11}), which will be discussed in Section~\ref{sec:conformal}.

We therefore adopt the following two‑step procedure named {\it modeling-then-calibration}:
\smallskip

(I)~\textbf{Modeling}: obtain a working model $\pi$ from a training data set
$\data_{\train}$ (or use a pre‑existing model).

(II)~\textbf{Calibration}: select $\lambda$ based on a separate calibration data
set $\data_{\cali}$ to minimize an empirical version of $\risk(\lambda)$
(defined precisely in equation \eqref{eq:asy-loss} below).
\smallskip

For the sake of simplicity, we assume that \(\bm{x}\) in $\data_\cali$ is uniformly distributed on the cube \([0,1]^d\). 
In this case, the indicators \(\mathbb{I}\bigl(\bm{x} \in \hat{\mathcal{X}}_R(\lambda),\; \bm{x} \notin \mathcal{X}_R\bigr)\) and \(\mathbb{I}\bigl(\bm{x} \not\in \hat{\mathcal{X}}_R(\lambda), \bm{x} \in {\mathcal{X}}_R\bigr)\) that appear in the risk translate directly into the Lebesgue measure (volume) of the symmetric difference \(\hat{\mathcal{X}}_R(\lambda) \Delta \mathcal{X}_R\) of the estimated set \(\hat{\mathcal{X}}_R(\lambda)\) and the true reliability set \(\mathcal{X}_R\).

Returning to the risk $\risk(\hat{\mathcal{X}}_R(\lambda))$ in \eqref{eq:risk11}, the critical difficulty in implementing calibration is that the events \(\{\bm x\in \mathcal{X}_R\}\) cannot be observed directly.  
To address this problem, we introduce an observable asymmetric loss function defined by 
\begin{equation}\label{eq:asy-loss-1}
\ell (\lambda, \bm x_i,Y_i ):= \mathbb{I}(\pi(\bm x_i)\ge 1-\lambda,\; Y_i=0) + (1-p_R)\,\mathbb{I}(\pi(\bm x_i) < 1-\lambda).
\end{equation}
The corresponding empirical loss for the calibration set \(\data_\cali\) is then  given by 
\begin{equation}\label{eq:asy-loss}
    {{\loss}_{n_\cali}(\lambda)} := \frac{1}{n_{\cali}}\sum_{i=1}^{n_\cali}\ell (\lambda, \bm x_i,Y_i ). 
\end{equation}

Our first result shows that the minimizer of $\E\{ {\loss}_{n_\cali}(\lambda) \}$ always yields a risk defined in \eqref{eq:risk11} no larger than the classical approach's.

\begin{theorem}\label{lem:consistency}
    Let \(\lambda^*\) be the minimizer of the population loss {\(
    \E\{ \loss_{n_\cali}(\lambda)\} \)} with a pre-trained model $\pi$. 
 It holds that
    \[\risk(\lambda^*)\le\risk(1-p_R).\]
    Moreover, the Lebesgue measure of symmetric difference $\hat{\mathcal{X}}_R(\lambda^*) \Delta \mathcal{X}_R $ between the  set estimator \(\hat{\mathcal{X}}_R(\lambda^*) \) and  the reliability set  \(\mathcal{X}_R\) 
    is zero if one of the following conditions holds.
        \begin{enumerate}\spacingset{1.1}
        \item[(a)]  \(\pi(\bm{x}) = \Pr(Y=1 \mid \bm{x})\) almost surely.
        \item[(b)]  There exists a strictly increasing function \(h(\cdot)\) such that \(\pi(\bm{x}) = h(\Pr(Y=1 \mid \bm{x}))\).  
        \item[(c)]  The  estimator \(\pi \) separates \(\mathcal{X}_R\) and \(\mathcal{X}_R^c\) in the sense that \(\inf_{\bm{x} \in \mathcal{X}_R} \pi(\bm{x}) > \sup_{\bm{x} \in \mathcal{X}_R^c} \pi(\bm{x})\).
    \end{enumerate}
\end{theorem}
 Theorem~\ref{lem:consistency} indicates that by minimizing the population risk $\E\{\loss_{n_\cali}(\lambda)\}$ with respect to the threshold $\lambda$, one can effectively
reduce the generalization error $\risk(\lambda)$ compared with the plug-in method.
Moreover, it can recover the reliability set $\mathcal{X}_R$ under mild conditions that extend the model assumption from correct specification to useful misspecification.   
More specifically, case (a) assumes that the model is correctly specified, which is aligned with the assumption underlying the classical plug-in estimator.  
Case  (b) relaxes the requirement to order preservation, which is widely adopted in practice for perceptron learning \citep{kalai2009isotron}. This is naturally satisfied by models trained via rank-based objectives such as AUC maximization \citep{herschtal2004optimising} or RankNet \citep{burges2005learning}. For a univariate linear classifier (e.g., logistic regression), this condition is equivalent to sign consistency of the working model’s slope.  
Case (c) further ensures identifiability through a separation condition between reliable and unreliable units in the output space of $\pi$, emphasizing that a perfect probabilistic model is unnecessary. Instead, it suffices to separate \(\bm x\) in \(\mathcal{X}_R\) and \(\mathcal{X}_R^c\).

Let $\lambdaCali$ minimize the empirical loss in \eqref{eq:asy-loss}, then the next result shows that the Lebesgue measure $\mu$ of the symmetric difference   $\hat{\mathcal{X}}_R(\lambdaCali) \Delta \hat{\mathcal{X}}_R(\lambda^*)$  of the sets  $\hat{\mathcal{X}}_R(\lambdaCali)$ and  $\hat{\mathcal{X}}_R(\lambda^*)$  converges  to zero with rate $O_P(n_{\cali}^{-1})$.  
Before formally stating the result, we introduce the following three regularity conditions.

\begin{assumption}\label{ass:11}
    There exists a constant $c_0 > 0$ such that $|\Pr(Y=1\mid\bm x) - p_R| \ge c_0$ almost everywhere (a.e.) on  $[0,1]^d$.
\end{assumption}

\begin{assumption}\label{ass:51} 
$\Pr(Y=1\mid\bm x)=g(\s(\bm x))$ for some measurable function $g$, where $\s (\bm x ) = 1 - \pi   (\bm x ) $. There exists a unique interval $[c^*_{-},c^*_{+}]\subseteq [0,1]$ such that 
\begin{equation}\label{eq:monoto}
  (p_R- \E\{\Pr(Y=1\mid\bm x)\mid \s(\bm x)=u\})(u - c^*_+) > 0 
\end{equation}
almost everywhere on  $[0,1]\setminus[c^*_{-},c^*_{+}]$ (if $c_{-}^*=c^*_+$, the interval degenerates to a single point).
\end{assumption}

\begin{assumption}\label{ass:21}
   Let $\bm Z$ be a random variable with a uniform distribution on $[0,1]^d$. The score variable $\s(\bm Z) := 1 - \pi(\bm Z)$ has  a continuous density $f_{\s}$ on $[0,c_-^*)$ and $(c_+^*,1]$ and  there exists a positive constant $c$ such that $c^{-1}\le f_{\s}(u)\le c$ for a.e.  $u\in [0,1]\setminus [c^*_{-},c^*_{+}]$ and $f_\s(u)=0$ for a.e. $u\in [c^*_{-},c^*_{+}]$.
\end{assumption}

 Condition~\ref{ass:11} is also known as the Massart condition \citep{massart2006risk} and 
 corresponds to the low noise condition introduced in \cite{tsybakov2004optimal}. 
 This regularity assumption ensures that there is no accumulation of probability mass near the decision boundary. 
 Condition~\ref{ass:51} is a standard and widely accepted regularity assumption in statistical learning \citep[see, for example,][]{zadrozny2002transforming,nalbantov2019note}.
It is often referred to as sign-consistency or classification calibration, and it is theoretically guaranteed when the estimator is trained by minimizing a strictly proper scoring rule, such as the cross-entropy or logistic loss \citep{gneiting2007strictly}. 
This condition is used to ensure that the risk function {$\lambda \mapsto \E\big \{ \loss_{n_\cali}(\lambda)\big \} $} is convex such that its minimizer  $\lambda^*$ is well defined. 
It automatically holds when the working model is monotonically related to $\Pr(Y=1\mid\bm x)$ via an unknown relationship.
Condition~\ref{ass:21} is a regularity assumption on the distribution of the score. It rules out discrete score distributions and guarantees that the score density is bounded and bounded away from zero outside the interval of optimal thresholds. When the optimal threshold is not unique, the condition allows for a gap in the score distribution on $[c_-^*,c_+^*]$.
 Conditions~\ref{ass:51} and~\ref{ass:21}  are formulated in a way to accommodate cases where the threshold is not unique. This situation appears, for example,  under the separation scenario described in Theorem~\ref{lem:consistency}(c).

\begin{theorem}\label{thm:calibration}
    Let $\pi$ be a fixed working model (or conditioned on $\data_\train$). 
    Suppose that Conditions~\ref{ass:11}--\ref{ass:21} hold, and further assume that the data $\data_{\cali}$ are i.i.d..  The volume (with respect to the Lebesgue  measure $\mu$)  of the symmetric difference $\hat{\mathcal{X}}_R(\lambdaCali) \Delta \hat{\mathcal{X}}_R(\lambda^*)$, satisfies:
    \begin{equation}     \mu\big(\hat{\mathcal{X}}_R(\lambdaCali) \Delta \hat{\mathcal{X}}_R(\lambda^*)\big) = O_P\big( n_{\cali}^{-1} \big).
    \end{equation}
\end{theorem}

Theorem~\ref{thm:calibration} establishes that the calibrated set estimator converges to the optimal calibrated set at a rate of $O_P(n_{\cali}^{-1})$. 
This is a parametric rate even if the working model is estimated via nonparametric procedure. The intuition is that calibration operates on a single scalar parameter $\lambda$ rather than on the full $d$-dimensional surface $\Pr(Y\mid\bm x)$. 
Once the working model is fixed, the remaining problem is therefore one-dimensional, and, under the conditions of Theorem~\ref{thm:calibration}, the calibration error translates into an 
$O_P(n_{\cali}^{-1})$ error in the estimated set.

A practical consequence of this rate is that calibration is a low-cost safeguard. Even a modestly sized calibration set can substantially improve the decision quality of a pre-trained or historically available working model. 
In settings where a working model already exists (for example, derived from a previous study, engineering heuristics, or a related task), calibration provides a principled way to adapt it to the current decision problem at minimal additional data cost.

We close this section with a brief discussion of the classical setting where the working model is correctly specified. Although our framework is designed for the misspecified regime, it is natural to ask whether calibration incurs an
unnecessary cost when the model happens to be correct. 
In the Supplement we show that the cost is indeed negligible. 
Under regularity conditions,  both the calibrated estimator with the
split \(n_{\cali} \asymp n_{\train}^{2/(2+W)}\log n_{\train}\) and the classical plug‑in approach achieve the same convergence rate of order \(O_P(n_{\train}^{-2/(2+W)})\), where \(W \in (0,2)\) is the bracketing entropy exponent of the model class 
(for a fair comparison, we use \(n_{\train} + n_{\cali}\) sample points for the classical approach).
Hence, even under correct specification, reserving a fraction of the data for calibration comes at no first-order asymptotic cost. 
Full statements and proofs are given in Section S.9 of the Supplement.

\section{Conformal Inference for Reliability Set Estimation}\label{sec:conformal}

Section~\ref{sec:method} shows that calibration reduces the total decision risk. 
In safety-critical applications, however, a small total risk alone may not be 
sufficient to ensure a trustworthy decision. Since the quality of the working 
model $\pi$ is generally unknown, a set estimator with a small total risk may 
still incur an unacceptably large false inclusion risk, that is, the risk of 
certifying unsafe inputs as reliable.
We formalize this as the type~I risk of the decision rule in the same spirit of hypothesis testing:
    \begin{align}
    \label{eq:risk-decomp1}
       \risk_1(\lambda) &= \E\bigl\{\bigl(p_R - Y\bigr) \, \mathbb{I}\bigl(\bm{x} \in \hat{\mathcal{X}}_R(\lambda), \bm{x} \notin {\mathcal{X}}_R\bigr) \bigr\}. 
    \end{align}
Let \((\bm{x}_{n_\cali+1}, Y_{n_\cali+1})\) be a new data point independent of \(\data_{\cali}\) and identically distributed.
Our goal in this section is to find a threshold \(\lambda\), computed from \(\data_\cali\), with the risk control guarantee:
\begin{equation}\label{eq:risk}
    \E\{\risk_1(\lambda)\} 
    := \E\Bigl\{\bigl(p_R - Y_{n_\cali+1}\bigr) \, \mathbb{I}\bigl(\bm{x}_{n_\cali+1} \in \hat{\mathcal{X}}_R(\lambda), \bm{x}_{n_\cali+1} \notin {\mathcal{X}}_R\bigr) \Bigr\} 
    \le \alpha,
\end{equation}
where \(\alpha > 0\) is a user-specified error rate.

The conformal risk control (CRC) framework, introduced by \citet{angelopoulos2022conformal} and building on the broader conformal prediction literature~\citep{vovk2005algorithmic}, provides a natural tool for this task. Rather than minimizing a combined loss, CRC uses a held-out calibration set to select a threshold that directly bounds the expected loss of a decision rule at a user-specified level \(\alpha\), with finite-sample validity under exchangeability. Subsequent work has extended the framework to more general loss structures and multi-step decision problems \citep{angelopoulos2023conformal}.
However, direct application of CRC to reliability set estimation fails, as the direct loss in \eqref{eq:risk} contains the true set membership \(\mathbb{I}(\bm{x} \in \hat{\mathcal{X}}_R(\lambda), \bm{x} \notin \mathcal{X}_R)\) which is unobservable. 
To make CRC applicable, we require an observable surrogate loss that (i) can be computed on the calibration data and (ii) provides a valid upper bound for the unobservable error of interest. 
The following lemma provides such an observable surrogate.

\begin{lemma}\label{thm:loss_upper}
Define the generalized risk (for a measurable function $g$) by 
\begin{equation}\label{hol5}
    \tilde{\risk}(g) := \E\bigl\{\mathbb{I}(g(\bm{x}) \ge 1, Y=0) + (1-p_R),\mathbb{I}(g(\bm{x}) < 1)\bigr\}.
\end{equation} 
Then for any \(\lambda\), the type~I risk in \eqref{eq:risk-decomp1} satisfies
\begin{equation}\label{eq:general_bound}
    \risk_1(\lambda) \le \E\{\loss_{n_{\cali}}(\lambda)\} - \min_g \tilde{\risk}(g).
\end{equation}
\end{lemma}

With this surrogate bound in place, we can now refine a CRC procedure for reliability set estimation.
Note that $0\le \min_g\tilde \risk(g)\le \min_{\lambda}\E(\loss_{n_\cali}(\lambda))$, and  that it follows from the proof of Theorem~\ref{lem:consistency} in the Supplement that the bound is sharp. 
To ensure strict control of the risk in  \eqref{eq:risk}, we propose to select the threshold, denoted by $\lambdaCRC$,   as
\begin{equation}\label{eq:hat-lambda}
   \lambdaCRC := \sup \Big\{ \lambda \in [c_+^*, 1] ~~\Big |~~ \frac{1}{n_\cali+1}\sum_{i=1}^{n_\cali}\ell^{\uparrow}(\lambda,\bm x_i,Y_i) + \frac{2}{n_\cali+1}  \le \alpha \Big\}, 
\end{equation}  
where $\ell^{\uparrow}(\lambda,\bm x_i,Y_i)=\sup_{c_+^*\le \gamma\le\lambda}\ell(\gamma,\bm x_i,Y_i)$
 is the monotone refinement of the loss function and $c_+^*$ is the constant from  Condition~\ref{ass:21}.
 Here, monotonicity is required to translate a pointwise evaluation into a uniform guarantee for the data-driven choice of $\lambda$. 
 
 Note that the optimal solution of $\min_g\tilde \risk(g)$ is the Bayes estimator $g^*(\bm{x}) = \mathbb{I}\{\Pr(Y=1\mid\bm{x}) > p_R\}$. 
 The level $\alpha$ in \eqref{eq:hat-lambda} cannot be less than $\min_g\tilde \risk(g)$, otherwise the resultant reliability set will be (asymptotically as $n_{\rm cal} \to \infty$) empty.   
 We use $c_+^*$ in the definition of the monotone refinement and of $\hat \lambda^\diamond$, since it is a minimizer of $\E\{\loss_{n_\cali}(\lambda)\}$ and this choice also provides a clear presentation of our results. A different threshold can be used here as well, and the subsequent result remains valid after updating the constant accordingly. However, a poorly chosen bound, for instance, one that is too small, may inflate the adjusted loss and render the reliability set empty even when \(\alpha\) exceeds the minimal risk. 
A fully data-adaptive version, which does not require the knowledge of the constant $c_+^*$, is given in \eqref{eq:hat-lambda1} below.

\begin{theorem}\label{thm:validity}
{Let Condition~\ref{ass:21} be satisfied.}
Assume  
\begin{align} \label{condalpha}
\Pr \Big (\frac{1}{n_\cali+1}\sum_{i=1}^{n_\cali}\ell^{\uparrow}(c_+^*,\bm x_i,Y_i)\le \alpha-\frac{3}{n_\cali+1} \Big )=1
\end{align}
and that the data $\data_{\cali}\cup\{(x_{n_\cali+1}, Y_{n_\cali+1})\}$ are i.i.d..
Given a pre-trained model $\pi$, the risk $\risk_1$ on the test set $(x_{n_\cali+1}, Y_{n_\cali+1})$ satisfies
     \begin{equation}
          \E\big\{\risk_1(\lambdaCRC)\big\}\le \E\Big \{\sup_{c_+^*\le\gamma\le \lambdaCRC}\ell(\gamma, \bm x_{n_\cali+1}, Y_{n_\cali+1})\Big \}\le \alpha.
     \end{equation}
\end{theorem}

We highlight that the risk control in Theorem~\ref{thm:validity} holds uniformly, independent of the underlying model quality and the size of the training and calibration sample.  Here,  condition \eqref{condalpha} is only used to ensure that $\alpha$ is a reasonable risk level.
Thus, the data scientist can safely use any threshold value in the interval  $[c_+^*, \lambdaCRC]$.
It is worth emphasizing that the monotone refinement typically does not introduce excessive conservatism, as the original loss is often near-monotone in engineering applications. 
Consequently, the supremum correction results in only a modest increase in the loss.
Moreover, the interval guarantee provides operational flexibility: practitioners may deliberately choose a lower threshold to obtain tighter reliability sets or reduce downstream computation costs, secure in the knowledge that the supremum-based control still bounds the worst-case risk.
Such flexibility is especially useful when the data-driven choice \(\lambdaCRC\) is overly cautious for the task at hand.
The next result establishes a lower bound for $\E\{\sup_{c_+^*\le\gamma\le \lambdaCRC}\ell_{n_\cali+1}(\gamma, \bm x_{n_\cali+1}, Y_{n_\cali+1})\}$, ensuring that the procedure is not overly conservative to satisfy the upper bound.

\begin{theorem}\label{thm:lower_bound}
   Assume that the conditions of Theorem~\ref{thm:validity} are satisfied and further assume that 
   there exists a constant $c\in[c_+^*,1]$ such that 
  \begin{align}
\Pr \Big ( \frac{1}{n_\cali+1}\sum_{i=1}^{n_\cali}\sup_{c_+^*\le\gamma\le c}\ell(\gamma,\bm x_i,Y_i)>\alpha-\frac{2}{n_\cali+1} \Big ) =1 ,
      \label{condalpha1}
   \end{align}
 then  we have
   \begin{equation}
       \E\Big \{\sup_{c_+^*\le\gamma\le \lambdaCRC}\ell(\gamma, \bm x_{n_\cali+1}, Y_{n_\cali+1})\Big \} \ge \alpha  - {\frac{2}{{n_\cali}+1}}.
   \end{equation}
\end{theorem}

Clearly, the lower bound fails when $\alpha\ge 2$ (In this case, the risk control in Theorem~\ref{thm:validity} is trivial). 
Thus, we require the additional assumption on $\alpha$  in \eqref{condalpha1} to ensure the lower bound is meaningful.
Together with Theorem~\ref{thm:validity}, our target is controlled arbitrary close to \(\alpha\) (as $n_\cali\to \infty$), which rules out the possibility that  $\lambdaCRC$ is overly conservative.
In scenarios where strict control of type-I risk via the bound $\E \{\sup_{c_+^*\le\gamma\le \lambdaCRC}\ell(\gamma, \bm x_{n_\cali+1}, Y_{n_\cali+1})\}$ is unnecessary, an asymptotic or less conservative guarantee can suffice. We therefore refine the construction in \eqref{eq:hat-lambda}. 

More precisely, by further assuming Conditions~\ref{ass:51} and \ref{ass:21},  Lemma~\ref{thm:loss_upper} can be strengthened to 
\begin{equation}\label{eq:exact-risk-excess}
    \risk_1(\lambda) = \E\{\loss_{n_{\cali}}(\lambda)\} - \min_{\lambda} \E\{\loss_{n_{\cali}}(\lambda)\}
    \qquad \text{for all } \lambda \ge c_+^*.
\end{equation}
(we prove this statement  in the proof of Theorem~\ref{thm:validity1} below in the Supplement).
This identity motivates us to target the excess risk directly. Specifically, we aim to find a data-driven threshold \(\lambdaExcess\) such that
\begin{equation}\label{eq:excess1}
    \E\Bigl\{ \ell(\lambdaExcess, \bm x_{n_\cali+1}, Y_{n_\cali+1}) - \min_{\lambda}\mathcal{L}(\lambda) \Bigr\} \le \alpha .
\end{equation}
To avoid the excessive conservatism of the pointwise monotone majorization used in the finite-sample construction, we replace it by \(\sup_{\lambdaCali\le\gamma\le\lambda}\loss_{n_\cali}(\gamma)\), the supremum of the average raw loss.
The scaling \(n_{\cali}/(n_{\cali}+1)\) and the additive \(1/(n_{\cali}+1)\) provide the conformal correction for the unobserved test point, ensuring a valid bound with \(n_{\cali}+1\) observations.
For practical computation, we replace the unknown fixed lower bound \(c_+^*\) by the data-adaptive \(\lambdaCali\), and use \(\loss_{n_\cali}(\lambdaCali)\) as a finite-sample proxy for the minimal expected surrogate loss.
To be precise, the data-driven threshold is defined by  
\begin{equation}\label{eq:hat-lambda1}
    \lambdaExcess := \sup \Bigl\{ \lambda \in [\lambdaCali, 1] \;\Big|\;
    \frac{n_\cali}{n_\cali+1}\sup_{\lambdaCali\le\gamma\le \lambda}\loss_{n_\cali}(\gamma)
    + \frac{1}{n_\cali+1} - \loss_{n_\cali}(\lambdaCali) \le \alpha \Bigr\},
\end{equation}
where we also define \(\lambdaExcess = -1\) if the defining set in \eqref{eq:hat-lambda1} is empty.
 Theorem~\ref{thm:validity1} establishes its asymptotic validity of the data-driven threshold $  \lambdaExcess$.

\begin{theorem}\label{thm:validity1}
  Assume that Conditions~\ref{ass:11}--\ref{ass:21} are satisfied and further assume that the data  {$\data_{\cali}\cup\{(\bm x_{n_\cali+1},Y_{n_\cali+1})\}$} are i.i.d.. Then the risk of the threshold  $\lambdaExcess$ in  \eqref{eq:hat-lambda1} with a pre-trained model $\pi$ satisfies
     \begin{equation}
          \limsup_{n_\cali\to\infty}\E \{\risk_1(\lambdaExcess) \}=\limsup_{n_\cali\to\infty}\E\left\{\ell(\lambdaExcess,\bm x_{n_\cali+1}, Y_{n_\cali+1})-\min_{\lambda}\E(\loss_{n_\cali}(\lambda))\right\}\le \alpha.
     \end{equation}
\end{theorem}
Theorem~\ref{thm:validity1} establishes that our conformal calibration procedure controls the excess risk relative to the best attainable threshold within the working model, even under misspecification. From the proof of Theorem~\ref{thm:validity1}, one can further obtain a sharper description of the asymptotic risk when the target level is moderate. If \(\alpha + \min_\lambda\E\{\loss_{n_\cali}(\lambda)\} < \E\{\loss_{n_\cali}(1)\}\), then the limit exists and equals \(\alpha\) exactly, i.e.,
\(\lim_{n_{\cali}\to\infty} \E\{\risk_1(\lambdaExcess)\} = \alpha\),
showing that the proposed threshold achieves the prescribed risk level, while avoiding unnecessary conservatism.
This guarantee is precisely what renders the reliability set trustworthy in practice: regardless of whether the working model \(\pi(\bm{x})\) approximates \(\Pr(Y=1\mid\bm{x})\) well, the Type~I risk is bounded by \(\alpha\) asymptotically.

\section{Optimal Design for Reliability Region Estimation}\label{sec:design}

The preceding Sections~\ref{sec:method} and~\ref{sec:conformal} show that 
calibration can mitigate the effects of an imperfect working model by selecting 
a data-driven decision threshold. 
Its effectiveness, however, still depends on the quality of the working model $\pi$. When the model is correctly specified, classical optimal design theory~\citep{atkinson2007sas} provides efficient 
strategies for data collection. In the setting of insufficient model information considered here, such model-based designs may no longer be appropriate. 
Instead, a decision-oriented design is needed that allocates observations to regions that are most relevant for accurate estimation of the reliability set.

To motivate such a design, we highlight that the reliability set $\mathcal{X}_R$ is defined by
$\Pr(Y=1 \mid \bm{x}) \ge p_R$.
Whether a point is correctly included in or excluded from the estimated set
depends on how $\pi$ behaves near the decision boundary
$\{\bm{x}\mid \Pr(Y=1 \mid \bm{x}) = p_R\}$.
Far from this boundary, even a large estimation error in $\pi$ is unlikely to
change the decision; near the boundary, even a small error can be decisive.
This motivates concentrating observations where the decision is most sensitive. On the other hand, 
in  practice,  experimenters additionally require an understanding of the response
surface inside the reliability region for downstream analysis.
We therefore define the critical region
$\mathcal{S}_{\dense} = \mathcal{X}_R \cup \mathcal{T}_w(\partial\mathcal{X}_R)$,
a tubular neighborhood of the boundary enlarged to cover the reliable set,
and denote its complement $\mathcal{S}_{\sparse}$.
Our design principles are:

\smallskip
    (I) {\bf Space-filling on $\mathcal{S}_{\dense}$.} 
        A space-filling design provides accurate recovery of $\pi$ where the
        decision is most sensitive and also supports downstream modelling inside
        $\mathcal{X}_R$; it is minimax-optimal for smooth function estimation on a compact set  \citep{XieFang2000,BiedermannDette2001}. 
        
  (II) {\bf Controlled allocation to $\mathcal{S}_{\sparse}$.}
        The expected proportion of points in $\mathcal{S}_{\sparse}$ is bounded by
        a pre-specified level. In this region, misclassification is inherently
        unlikely, so we may limit the budget without compromising the reliability
        of the final decision.  
\smallskip

Implementing these two principles faces a fundamental circularity:
the critical region $\mathcal{S}_{\dense}$ depends on the unknown reliability
set and its boundary, and hence cannot be identified before data are collected.
We therefore partition the entire design space $\mathcal{X}$ into finitely many
grid cells and learn adaptively which cells should receive further observations.
During an initial warm-up stage, every cell is explored to obtain preliminary
information on its average response probability. Afterwards, we treat the cells
as arms in a Thompson sampling procedure \citep{sutton2018RL}: at each round $t$, cells whose posterior draws exceed a threshold $\tau$ form the candidate set $S(t)$ and continue to
receive observations, while cells with persistently low response probabilities
are gradually screened out.
Algorithm~\ref{alg:bandit_learning_warmup} formalizes this strategy in two stages.

\medskip

\begin{algorithm}[H]
\caption{Grid-Bandit Active Learning}
\label{alg:bandit_learning_warmup}
\spacingset{0.9}
\SetKwInOut{Input}{Input}
\SetKwInOut{Output}{Output}

\textbf{Input:}{ Total rounds $T_{\train}$, Grid division constant $K$, Bandit threshold $\tau$, Warm-up proportion $\nu$}\;
\textbf{Output:}{ Training data set $\data_{\train}$.}

\BlankLine
\textbf{Initialization:} {Define $T_{\rm warm} = \lfloor \nu T_\train \rfloor$,} partition $\mathcal{X}$ into $M = K^d$ disjoint grid cells $\{\region_1, \dots, \region_M\}$ and initialize the training set
$\data_{\train} \leftarrow \emptyset$,$\alpha_1(0)=\ldots=\alpha_m(0)=1, \beta_1(0)=\ldots=\beta_m(0)=1$\;

\BlankLine
{
\For{$t = 1$ \KwTo $T_{\rm warm}$}{
    \ForEach{$m \in \{1, \dots, M\}$}{
        {Add a design point $\bm x_m$ drawn from a uniform distribution on the cell $\region_{m}$\;
        Observe response $y_m \in \{0, 1\}$ and  update $\data_{\train} \leftarrow \data_{\train} \cup \{(\bm x_m, y_m)\}$\;}
       {Update priors: $\alpha_m(t) \leftarrow \alpha_m(t-1) + y_m$, $\beta_m(t) \leftarrow \beta_m(t-1)+1-y_m$\;}
    }
}
}

\BlankLine

\For{$t = T_{\rm warm} + 1$ \KwTo $T_{\train}$}{
    \For{$m = 1$ \KwTo $M$}{
        Sample $\tilde{\theta}_m(t) \sim \text{Beta}(\alpha_m(t-1), \beta_m(t-1))$\;
    }

    Identify candidate cells $S(t) \leftarrow \big \{m \in \{ 1, \ldots , M\} \mid \tilde{\theta}_m(t) \ge \tau \big \}$ (if $S(t) = \emptyset$ regenerate $\tilde{\theta}_m(t)$ for $m=1,\ldots,M$ again)\;

    \ForEach{$\breve m \in S(t)$}{
        Add a design point $\bm x$ drawn from a uniform distribution on the cell $\region_{\breve m}$\;
        Observe response $y \in \{0, 1\}$ and update the training set
        $\data_{\train} \leftarrow \data_{\train} \cup \{(\bm x, y)\}$\;
        
        \eIf{$y = 1$}{
            $\alpha_{\breve m}(t+1) \leftarrow \alpha_{\breve m}(t) + 1$\;
        }{
            $\beta_{\breve m}(t+1) \leftarrow \beta_{\breve m}(t) + 1$\;
        }
    }
}

\BlankLine
\Return{$\data_{\train}$}\;
\end{algorithm}
\bigskip

Clearly, Algorithm~\ref{alg:bandit_learning_warmup} implements the two principles above.
Within repeatedly selected cells, design points are spread over the cell,
yielding a space-filling design over the cells covering
$\mathcal{S}_{\dense}$ and hence supporting accurate estimation where it matters
most. In contrast, observations in $\mathcal{S}_{\sparse}$ arise mainly from
 the warm-up stage or occasional bandit mis-selections.
 To illustrate the operational dynamics of the proposed algorithm, we consider the following Example~\ref{eg:bandit-eff}.

\begin{figure}[h]
\centering\spacingset{0.9}
\includegraphics[width=\linewidth]{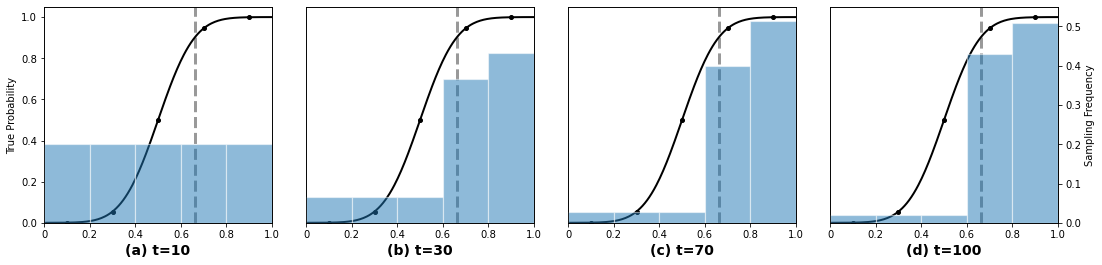}
\caption{\it Evolution of the design (displayed by the histograms) in the warm-up phase of Algorithm \ref{alg:bandit_learning_warmup} in the one-dimensional Probit model at rounds {\ $t=10, 30, 70, 100$}. The {solid curves}  represent the true probability $\Pr(Y=1\mid x)$.  The {vertical dashed line} demarcates the high-probability region $\mathcal{X}_R$.}
\label{fig:bandit_process}
\end{figure}

\begin{example}\label{eg:bandit-eff}
    Consider a one-dimensional Probit model with $\Pr(Y=1\mid x) = \Phi(8x-4)$ on the design space $\mathcal{X}=[0,1]$. 
    The hyper-parameters in Algorithm~\ref{alg:bandit_learning_warmup} are set as $T_\train=100$, $K=5$, $\tau=0.80$, and $\nu=0.1$.
    
    In Figures \ref{fig:bandit_process}(a)-(d), we illustrate the evolution of the allocation of the design at distinct rounds: $t=10$, $30$, $70$, $100$, respectively. 
    In these plots, the solid curve represents the underlying probability function \(\Pr(Y=1 \mid x)\), and the histogram depicts the approximate design (sampling frequency) within each region. 
    Crucially, the vertical dashed line marks the boundary of the high-probability region \(\mathcal{X}_R\) (where \(\Pr(Y=1 \mid x) \ge 0.9\)). The area to the right of this dashed line corresponds to our target region \(\mathcal{X}_R\).
    
    In Figure \ref{fig:bandit_process}(a), we show the initial stage of full exploration (warm-up) at $t=10$, where the sampling budget is uniformly distributed across the design space. 
    Figure \ref{fig:bandit_process}(b) captures the transition phase at $t=30$, where a balance between exploration and exploitation emerges; notably, the sampling frequency begins to increase in the grid cells located to the right of the dashed line.
    In Figures \ref{fig:bandit_process}(c) and (d) ($t=70$ and $t=100$), the algorithm demonstrates strong exploitation behavior. 
    The sampling budget concentrates heavily within the region \(\mathcal{X}_R\) (to the right of the vertical dashed line). 
    Thus the algorithm effectively ``screens out'' the sub-optimal cells on the left, 
    while focusing resources on the targeted level set.
\end{example}

To give this a precise operational meaning, let \(\mathcal{M}_{\sparse}\) be the maximal subset of cells that are fully contained in the sparse region \(\mathcal{S}_{\sparse}\), and let \(\mathcal{M}_{\dense} = \mathcal{M} \setminus \mathcal{M}_{\sparse}\) denote the remaining cells covering the critical region. 
 We require a separation condition on the conditional probabilities of these cells, which is formalized in the following assumption.
\begin{assumption}\label{ass:81}
$\Delta : = \min_{m\in\mathcal{M}_{\dense}} \Pr(Y=1 \mid \bm{x} \in \region_m) - \max_{m\in\mathcal{M}_{\sparse}} \Pr(Y=1 \mid \bm{x} \in \region_m) > 0.$
\end{assumption} 
This condition is mild and can be satisfied when
the reliability set \(\mathcal{X}_R\) is convex and that
\(\|\nabla \Pr(Y=1 \mid \bm{x})\|\) is bounded away from zero on the complement
\(\mathcal{X}_R^c\).
Under these assumptions, the response probability decreases at a positive minimal
rate as one moves from the boundary into the sparse region.
Since the sparse region \(\mathcal{S}_{\sparse}\) is separated from the critical
region \(\mathcal{S}_{\dense}\) by a tubular neighborhood of width \(w>0\), there
exists a constant \(\delta>0\) such that
\(
\inf_{\bm{x} \in \mathcal{S}_{\dense}} \Pr(Y=1 \mid \bm{x}) \;-\;
\sup_{\bm{x} \in \mathcal{S}_{\sparse}} \Pr(Y=1 \mid \bm{x}) \;\ge\; \delta.
\)

Let $\mathscr{A}(t) = \{ S(t) \ne \mathcal{M}_{\dense}\}$, where  $S(t)$ is  the set of grid cells  selected by Algorithm~\ref{alg:bandit_learning_warmup}  at round $t$. 
Note that ${\mathscr A}^c(t)$ denotes the event that, at time $t$, Algorithm~\ref{alg:bandit_learning_warmup} correctly identifies all cells in ${\cal M}_\dense$. In contrast, on the event ${\mathscr A}(t)$, there may exist cells in $S(t) \cap {\cal M}_\dense$ that are not selected and cells in $S(t) \cap {\cal M}_\sparse$ that are selected. We refer to these incorrectly identified cells as \emph{mis-selected cells}.
If the proportion of rounds in which the algorithm fails to correctly identify
$\mathcal{M}_{\dense}$ is kept small, then Principle~II is satisfied by
construction and, because the algorithm allocates observations uniformly among
the selected cells, the design points concentrate on $\mathcal{M}_{\dense}$ in a
space‑filling manner, fulfilling Principle~I.
We now formalize this logic.

\begin{theorem}\label{thm:error_count_bound}
Assume that the partitioning $\region_1,\ldots,\region_M$ of Algorithm~\ref{alg:bandit_learning_warmup} satisfies Condition~\ref{ass:81}, set
\begin{align*}
\tau & =\frac 1 2\left(\min_{m\in\mathcal{M}_{\dense}} \Pr(Y=1 \mid \bm x \in \region_m)+\max_{m\in\mathcal{M}_{\sparse}} \Pr(Y=1\mid\bm x \in \region_m)\right) ,
 \end{align*}
  and choose
\(
\nu>
{(32\log T_{\train}/\Delta^2+1)}/{T_{\train}}
\). 
The expected proportion of rounds with mis-selected cells is bounded by:
\begin{equation}\label{eq:sparse_control}
    \E \Big\{\frac {1}{T_{\train}} \sum_{t=1}^{T_{\train}} \mathbb{I}(\mathscr{A}(t)) \Big \} \le {\nu} + {\frac{M}{T_\train} \Big ( 1 + \frac{16}{\Delta^2} \Big )}.
\end{equation}
\end{theorem}

Theorem~\ref{thm:error_count_bound} enables us to have a closer look at the effects of different factors, such as the number of grid cells $M$, and the design budget $T_\train$, on the expected allocation proportion in $\mathcal{M}_{\sparse}$.
Heuristically, the expected proportion of mis-selection is bounded by the exploration factor $\nu$ as the second term diminishes linearly as $T_{\train}$ converges to infinity.

At the end of this section, we briefly discuss the origin of the two terms in the bound~\eqref{eq:sparse_control}.
The first term, $\nu$, essentially corresponds to the exploration cost inherent in Thompson sampling.
Specifically, the algorithm must allocate a logarithmic number of samples to the  grid cells in $\mathcal{M}_{\sparse}$ to statistically verify that their average probabilities fall below the threshold $\tau$. 
The second term corresponds to the error associated with the distance between the Beta-distributed random variable $\tilde{\theta}_m$  generated  in Algorithm \ref{alg:bandit_learning_warmup} 
and the probability $\Pr(Y=1 \mid \bm x \in \region_m)$. 
As is common in the multi-armed bandit literature \citep[see, e.g.,][Theorem 1]{wang2018th}, the error can be decomposed into two sources: the randomness induced by the Beta distribution and the estimation error due to the finite numbers of successes $s_m$ and failures $f_m$ observed in the training data $\data_{\train}$ within cell  $\region_m$.
Consequently, we do not suggest a large value $\nu$ in practice as randomness from the beta distribution also leads to some kinds of exploration.

\section{Numerical Studies}\label{sec:numerical}

In this section, we evaluate the empirical performance of the proposed Grid-Bandit Active Learning strategy (Algorithm~\ref{alg:bandit_learning_warmup}) together with the proposed conformal inference framework in Section~\ref{sec:conformal}. We consider the estimation and inference problem for the reliability set with $p_R=0.6$ and $0.9$ on the design space $[0,1]^d$, respectively. 
The average performance and corresponding standard errors are calculated through 500 simulation runs.
All methods are implemented in Python 3.9.12 with default parameter settings.

\subsection{Synthetic Data Setups}

We consider two distinct data generating processes $\Pr(Y=1 \mid \bm{x})$ defined on the unit square $\mathcal{X} = [0, 1]^d$ for $d=2$  and $d=4$ (results for the latter are relegated to Section S.1 of the Supplement): 
\begin{itemize}
    \item \textbf{Truncated Logistic Regression (T-Logit)}. %
    The conditional probability is defined via a logistic link function:
$$\Pr (Y=1\mid \bm{x}) = \text{clip} \Big ( \big ({1 + \exp(- 15x_1-15x_2-1.5d+24.1)} \big )^{-1}, 0.35, 0.95 \Big ),$$
    where \(\text{clip}(v, a, b) = \max(a, \min(b, v))\) is a truncation operator. 
   This model represents a smooth but steep transition of probabilities. 

   \item  \textbf{Piecewise Constant (PWC)}. The conditional probability is defined via
$$
\Pr(Y=1\mid\bm{x}) = 0.2 + 0.2~  \mathbb{I}\big (x_1+x_2 \geq 0.7\big  ) + 0.55~  \mathbb{I}\big (x_1+x_2 \geq 1.5\big ).
$$
We use this model to test the methods under non-smooth conditions. 
\end{itemize}
We opt to use these parameter setups to ensure that the volumes of the reliability set $\mathcal{X}_R$ (with $p_R = 0.9$) are kept within approximately $10\%$–$15\%$. The working model $\pi$ is deliberately chosen to be a standard linear logistic
regression model. This choice reflects the practical reality that the experimenter rarely
knows the true response surface and must work with a tractable, interpretable
model that is at best a useful approximation.

To evaluate the effectiveness of both the proposed estimation procedure and the proposed design strategy, we compare our method against the classical plug‑in estimator. For each of these two estimation methods, we consider two design schemes for the training step: simple random sampling and the corresponding optimal design. Further details of these designs are given below.
\begin{itemize}
    \item \textbf{Simple Random Sampling (SRS).} As a classic baseline, we consider a standard approach where all $n$ samples $\bm{x}_i$ are drawn i.i.d. from a uniform distribution over the design space $[0,1]^d$.
    \item \textbf{Optimal Design (Opt).} When the plug-in method is adopted, we opt to use the optimal design for the working model proposed in \cite{dror2008sequential}. When the new method proposed in this paper is adopted, an optimal design is generated by Algorithm~\ref{alg:bandit_learning_warmup}.

\end{itemize}

The design budgets are controlled with total sample sizes \(n = 300\). 
For our proposed method, we split each budget into training and calibration sets in a 2:1 ratio. 
To strictly enforce the space-filling property within \(\mathcal{M}_{\dense}\), we employ a sequential minimax design~\citep{santner2003design} inside the selected cells rather than independent uniform draws in the warm-up stage of Algorithm~\ref{alg:bandit_learning_warmup}.
The warm-up proportion is $\nu = 0.2$, and the bandit threshold $\tau$ is set to $0.75$ and $0.55$ for $p_R = 0.9, 0.6$, respectively.
A further sensitivity analysis for these hyperparameters is provided in Figure~S.1 of the Supplement. 
These results indicate that our method is not very sensitive with respect to the choice of these hyperparameters.
For a fair comparison, classical methods use the full budget \(n\) for training and adopt the standard plug‑in estimator to construct the reliability set.

\subsection{Estimation Performance}\label{subsect:est-sim}

 Let $\hat{\mathcal{X}}_R$ denote the estimator of the reliable region obtained by different methods, where the working model $\pi$ is chosen as linear logistic regression. 
 The estimation performance is assessed using the volume of the symmetric difference (\textbf{S-Diff}),  $\mu(\hat{\mathcal{X}}_R \Delta \mathcal{X}_R)$, evaluated by $10000$ Monte Carlo points i.i.d. from a uniform distribution on the design space.

\begin{table}[h]
    \centering\spacingset{1.0}
    \caption{\it Empirical median and interquartile range (IQR, reported in parentheses) of the volume of the symmetric difference (S-diff) for different designs (SRS and Opt) and different analysis methods (classical Plug-in approach and our proposed method denoted by ``Ours''). The best performances are highlighted in boldface.}
    \label{tab:est_comp}
    \footnotesize 
    \setlength{\tabcolsep}{4.5pt} 
    \renewcommand{\arraystretch}{1.2} 

    \begin{tabular}{ c c cccc cccc}
    \toprule
     \multirow{2}{*}{$\bm{p_R}$} & \multirow{2}{*}{\textbf{}} & \multicolumn{4}{c}{\textbf{Plug-in}}& \multicolumn{4}{c}{\textbf{Ours}}\\
    \cmidrule(lr){3-6} \cmidrule(lr){7-10}
     & & \multicolumn{2}{c}{SRS} & \multicolumn{2}{c}{Opt} & \multicolumn{2}{c}{SRS} & \multicolumn{2}{c}{Opt} \\
    \midrule
    
     \multirow{2}{*}{0.6} 
    & T-logit & 0.0506 & (0.0505) & 0.1019 & (0.0686) & 0.0455 & (0.0434) & \textbf{0.0356} & (0.0343) \\ 
     & PWC & 0.0728 & (0.0616) & 0.1847 & (0.0597) & 0.0338 & (0.0412) & \textbf{0.0281} & (0.0349) \\ 
    \cmidrule(lr){2-10}
     \multirow{2}{*}{0.9} 
    & T-logit & 0.0999 & (0.0005) & 0.0997 & (0.0085) & 0.0393 & (0.0348) & \textbf{0.0344} & (0.0354) \\ 
     & PWC & 0.1250 & (0.0004) & 0.1165 & (0.0535) & 0.0352 & (0.0425) & \textbf{0.0292} & (0.0403) \\ 
    \bottomrule
    \end{tabular}
\end{table}

Table~\ref{tab:est_comp} reports the quantitative results for the working model \(\pi\) being moderately misspecified as a standard logistic regression. 
The proposed method consistently outperforms the classical plug‑in approaches, both in design and inference. 
Algorithm~\ref{alg:bandit_learning_warmup} yields substantially smaller volumes of the symmetric difference (S-diff)  across different data generating processes and target thresholds, confirming the theoretical results in Theorem~\ref{lem:consistency}. 
The improvement over the plug‑in estimator is significant.
For example, for the PWC with \(p_R=0.9\), S‑diff drops from 0.1250 (classical plug‑in with SRS) to 0.0353 (our method with SRS), and even below the plug‑in with optimal design (0.1165). 
Within our framework, the optimal design consistently beats SRS (e.g., 0.0292 vs. 0.0352), supporting the design rationale in Section~\ref{sec:design}.

\begin{figure}[h!]
    \centering\spacingset{1.0}
    \begin{subfigure}[b]{0.225\textwidth}
        \centering
        \includegraphics[width=\textwidth, height=\textwidth]{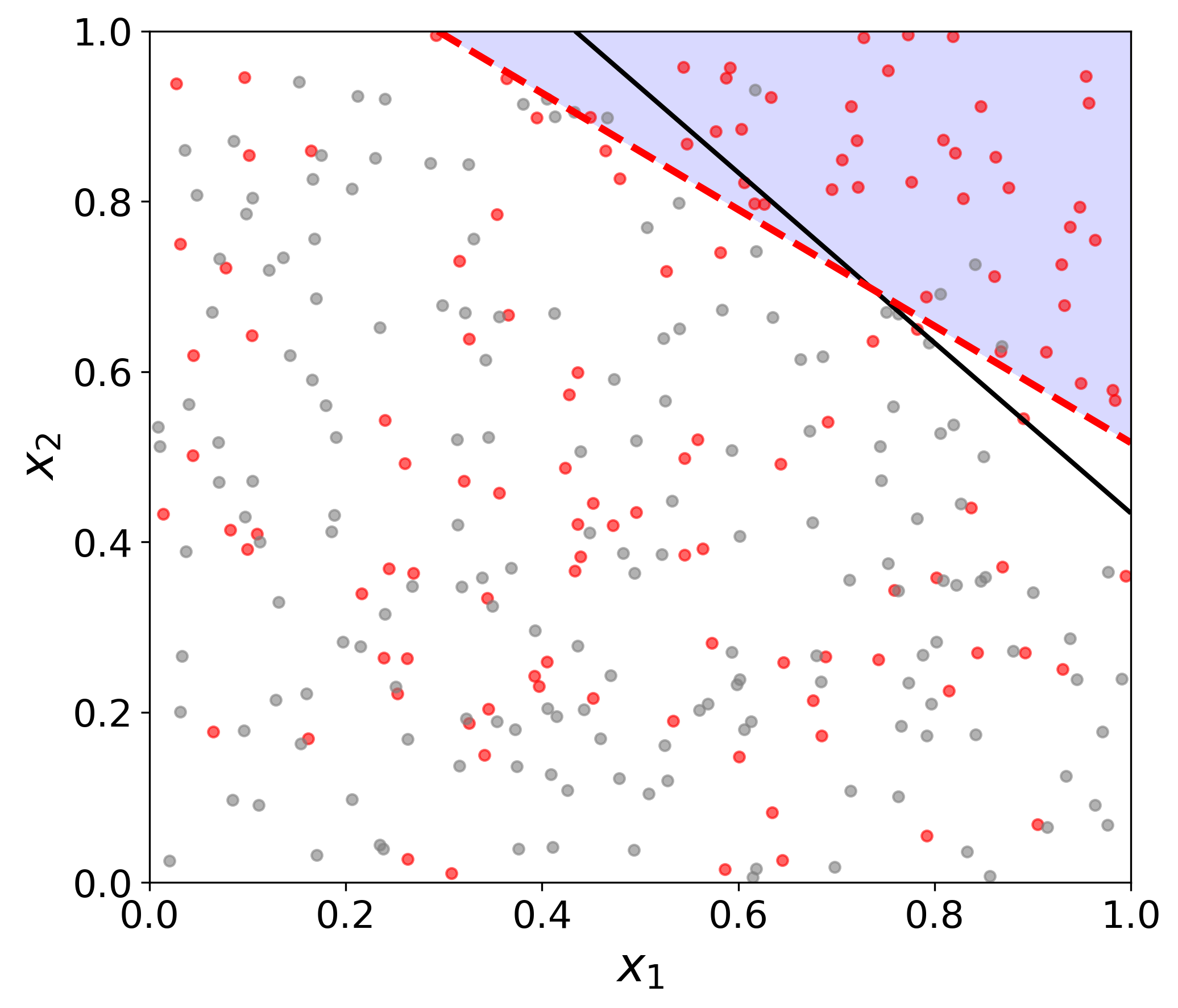} 
    \end{subfigure}
    \hfill 
    \begin{subfigure}[b]{0.225\textwidth}
        \centering
        \includegraphics[width=\textwidth, height=\textwidth]{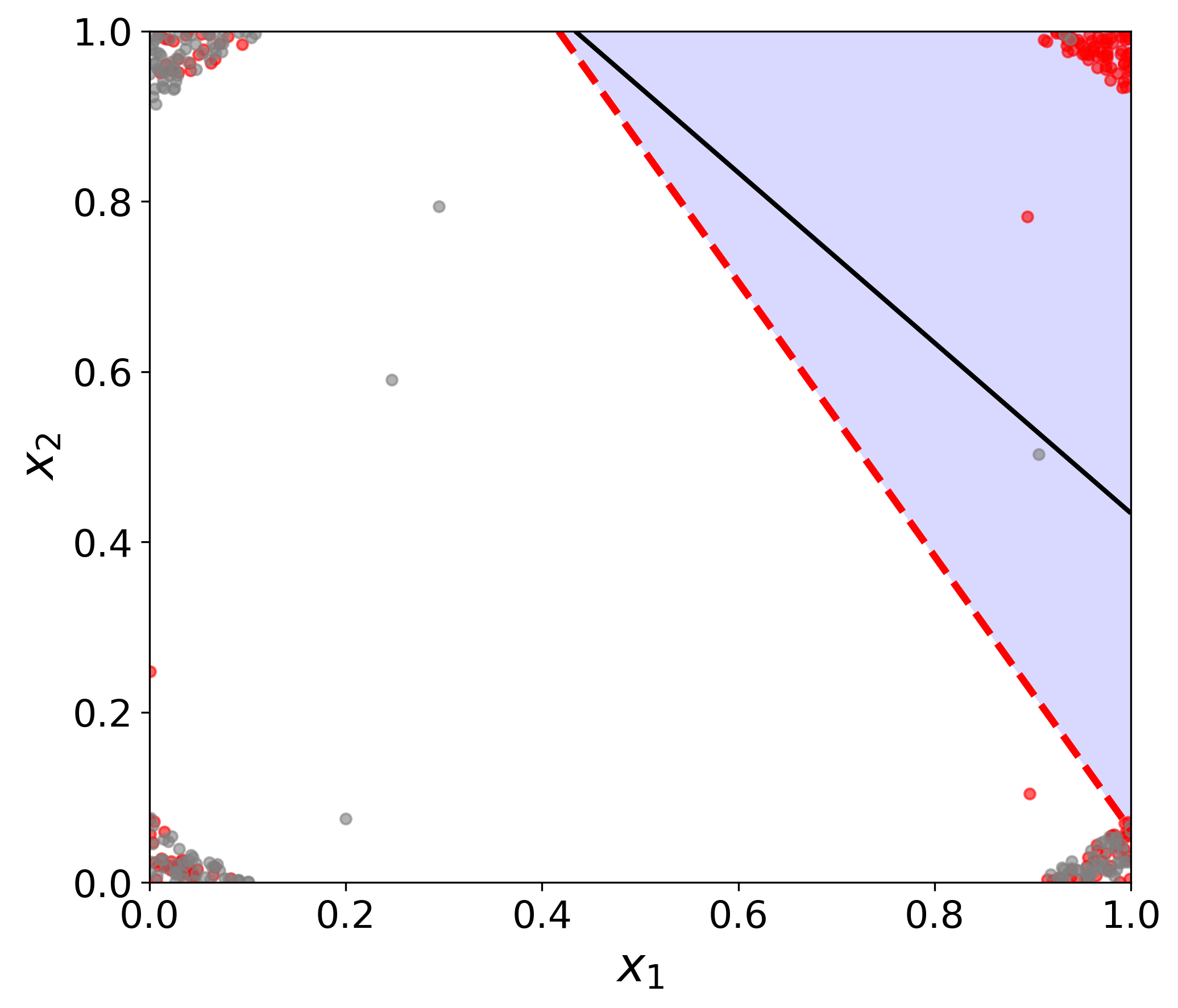}
    \end{subfigure}
    \hfill
    \begin{subfigure}[b]{0.225\textwidth}
        \centering
        \includegraphics[width=\textwidth, height=\textwidth]{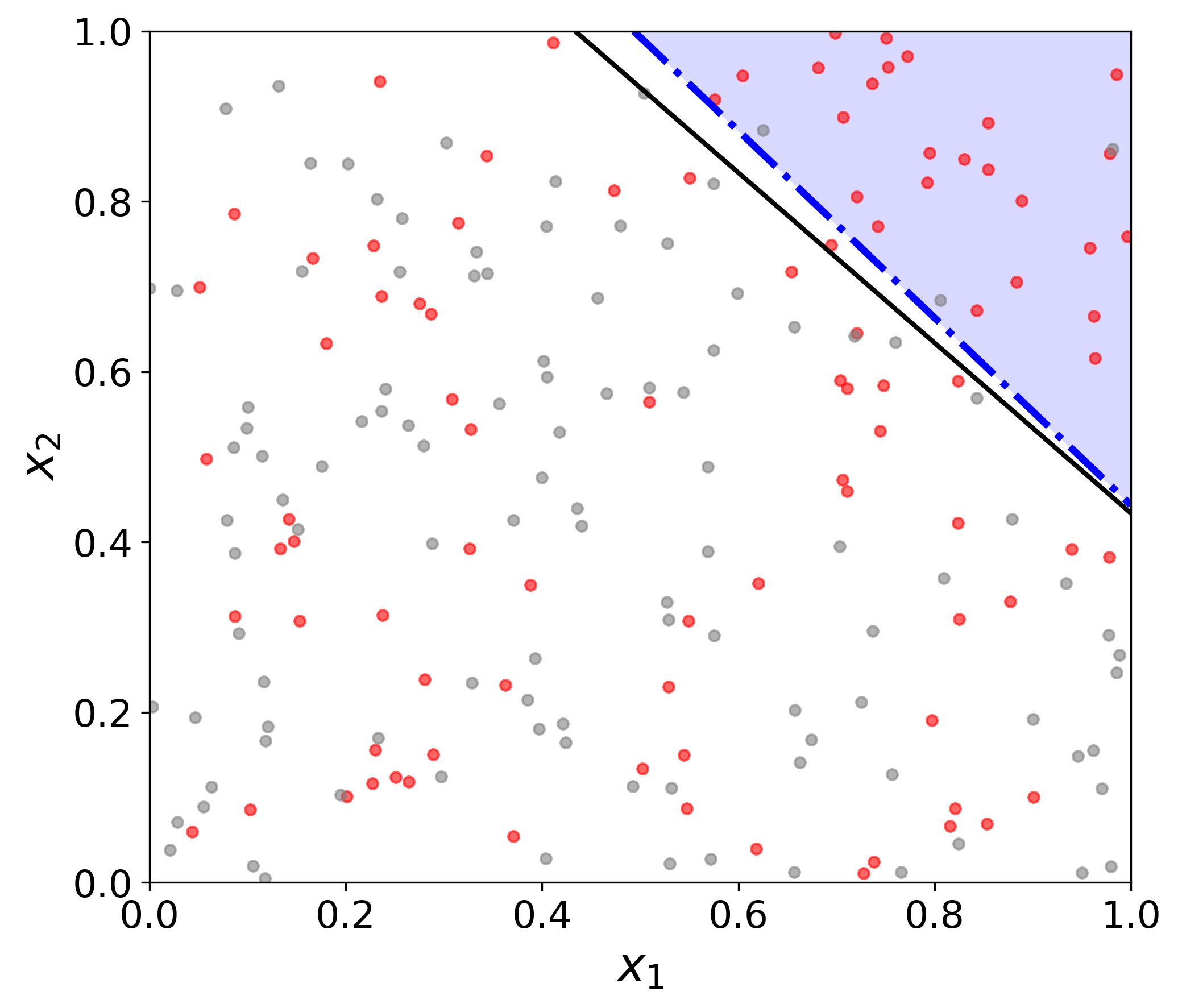}
    \end{subfigure}
    \hfill
    \begin{subfigure}[b]{0.225\textwidth}
        \centering
        \includegraphics[width=\textwidth, height=\textwidth]{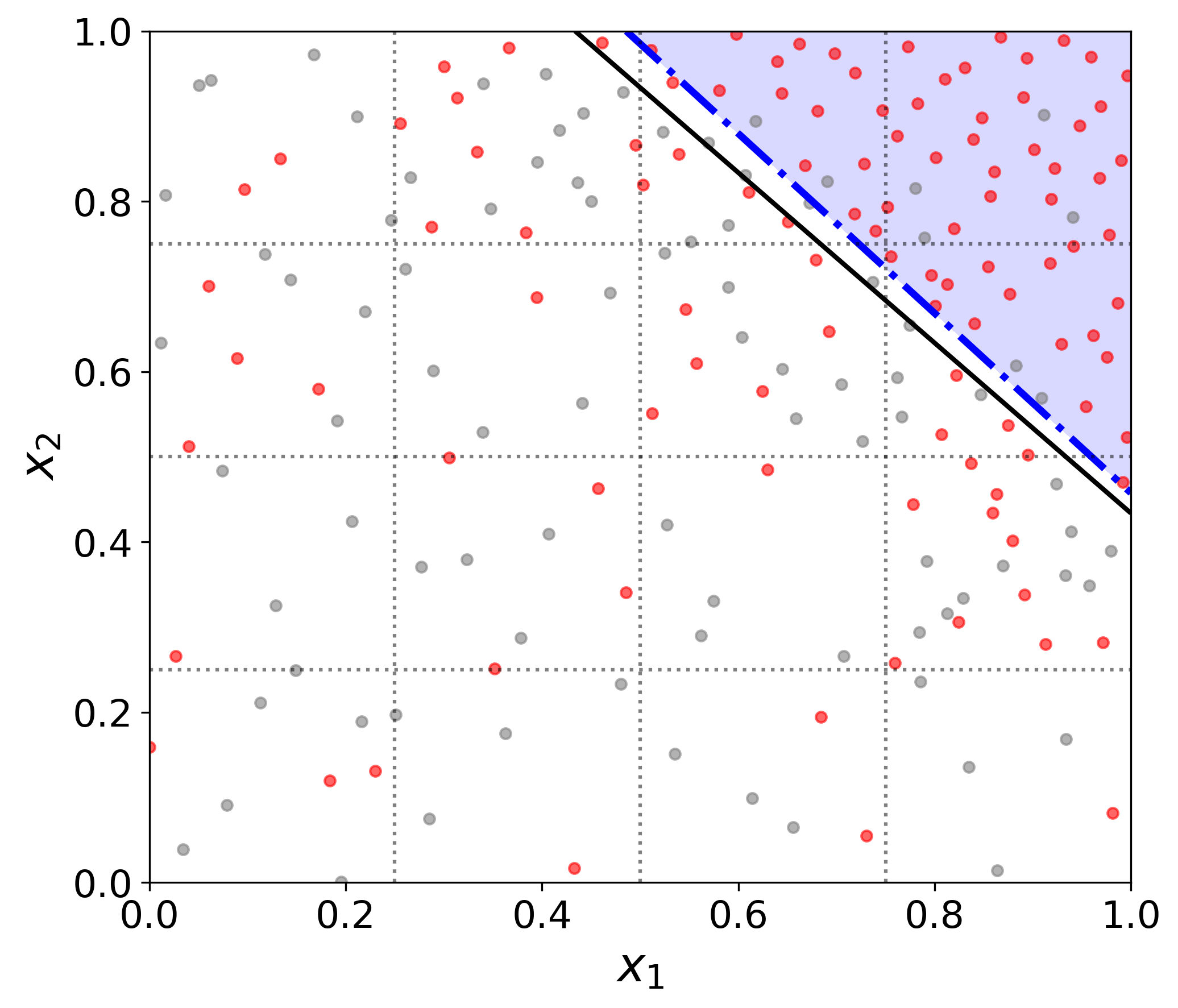}
    \end{subfigure}

    \vspace{1em}

    \begin{subfigure}[b]{0.225\textwidth}
        \centering
        \includegraphics[width=\textwidth, height=\textwidth]{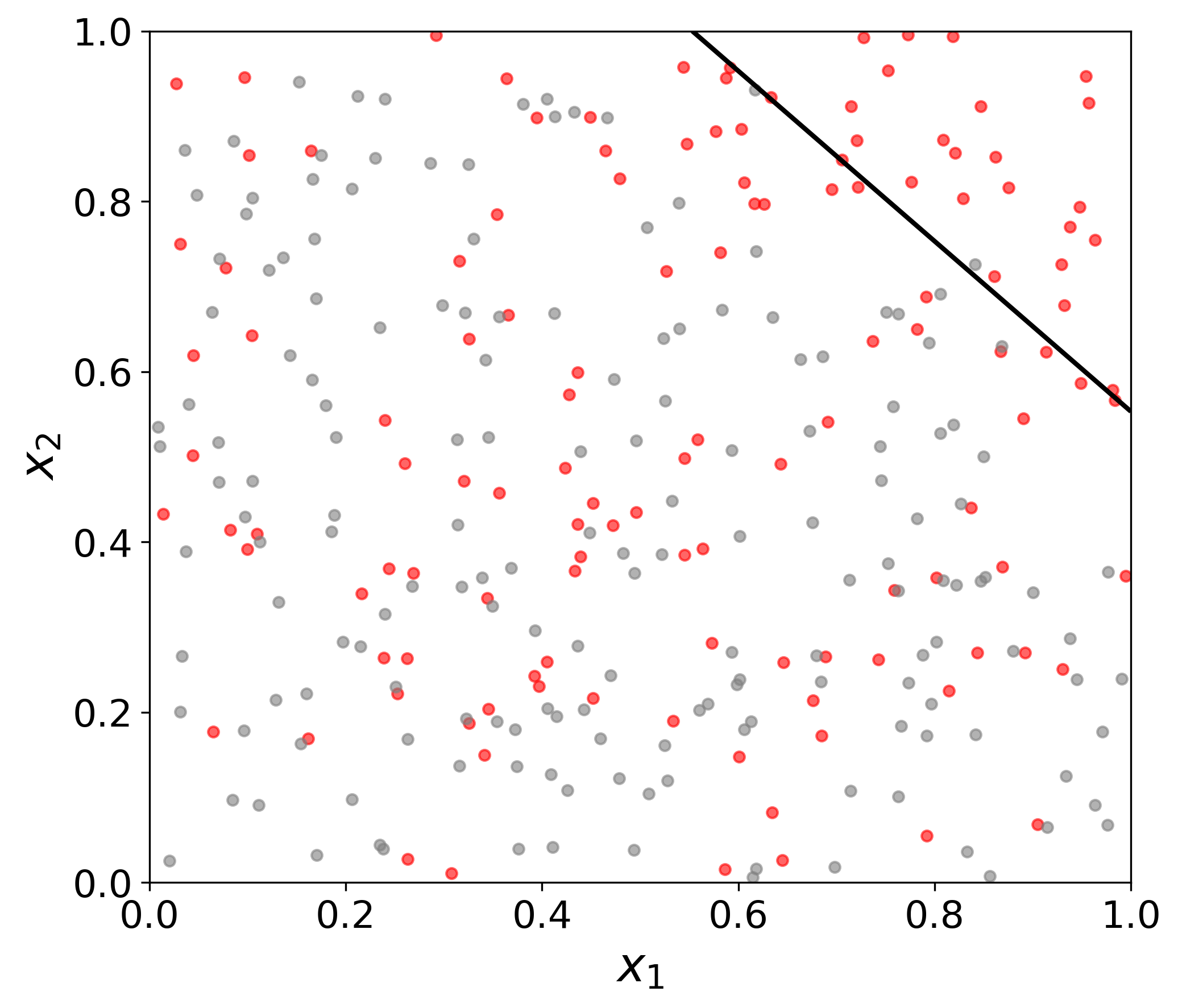}
    \end{subfigure}
    \hfill
    \begin{subfigure}[b]{0.225\textwidth}
        \centering
        \includegraphics[width=\textwidth, height=\textwidth]{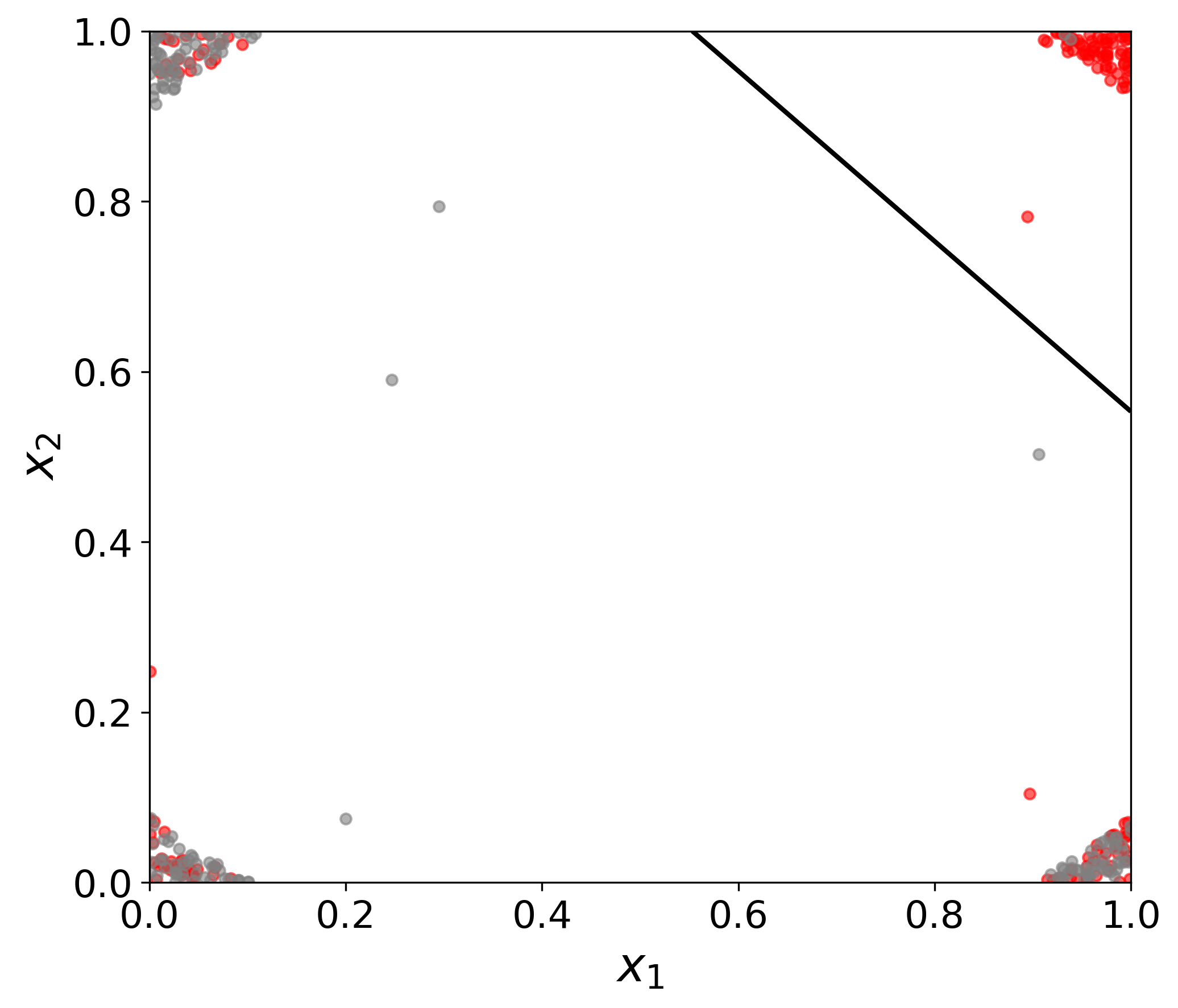}
    \end{subfigure}
    \hfill
    \begin{subfigure}[b]{0.225\textwidth}
        \centering
        \includegraphics[width=\textwidth, height=\textwidth]{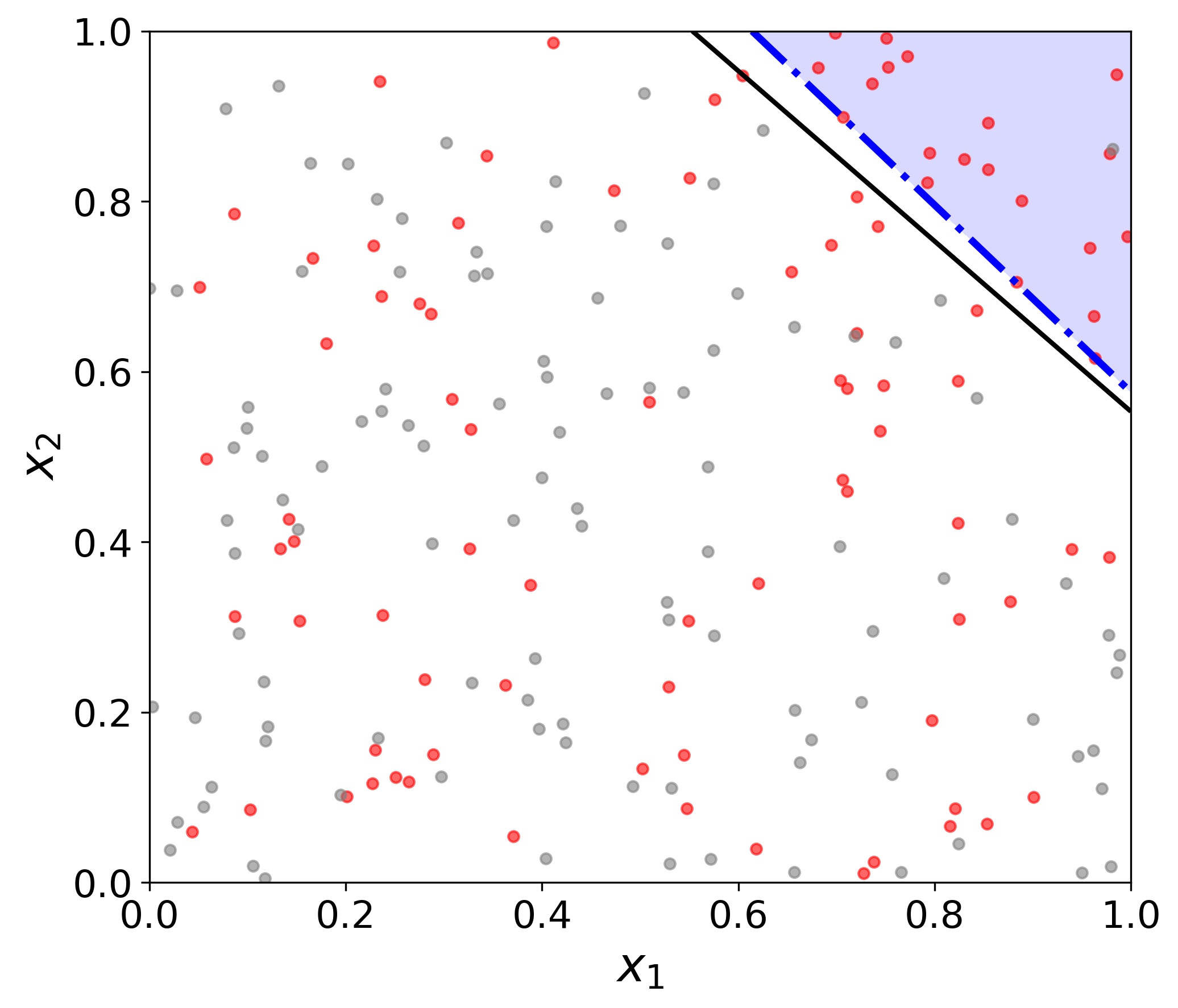}
    \end{subfigure}
    \hfill
    \begin{subfigure}[b]{0.225\textwidth}
        \centering
        \includegraphics[width=\textwidth, height=\textwidth]{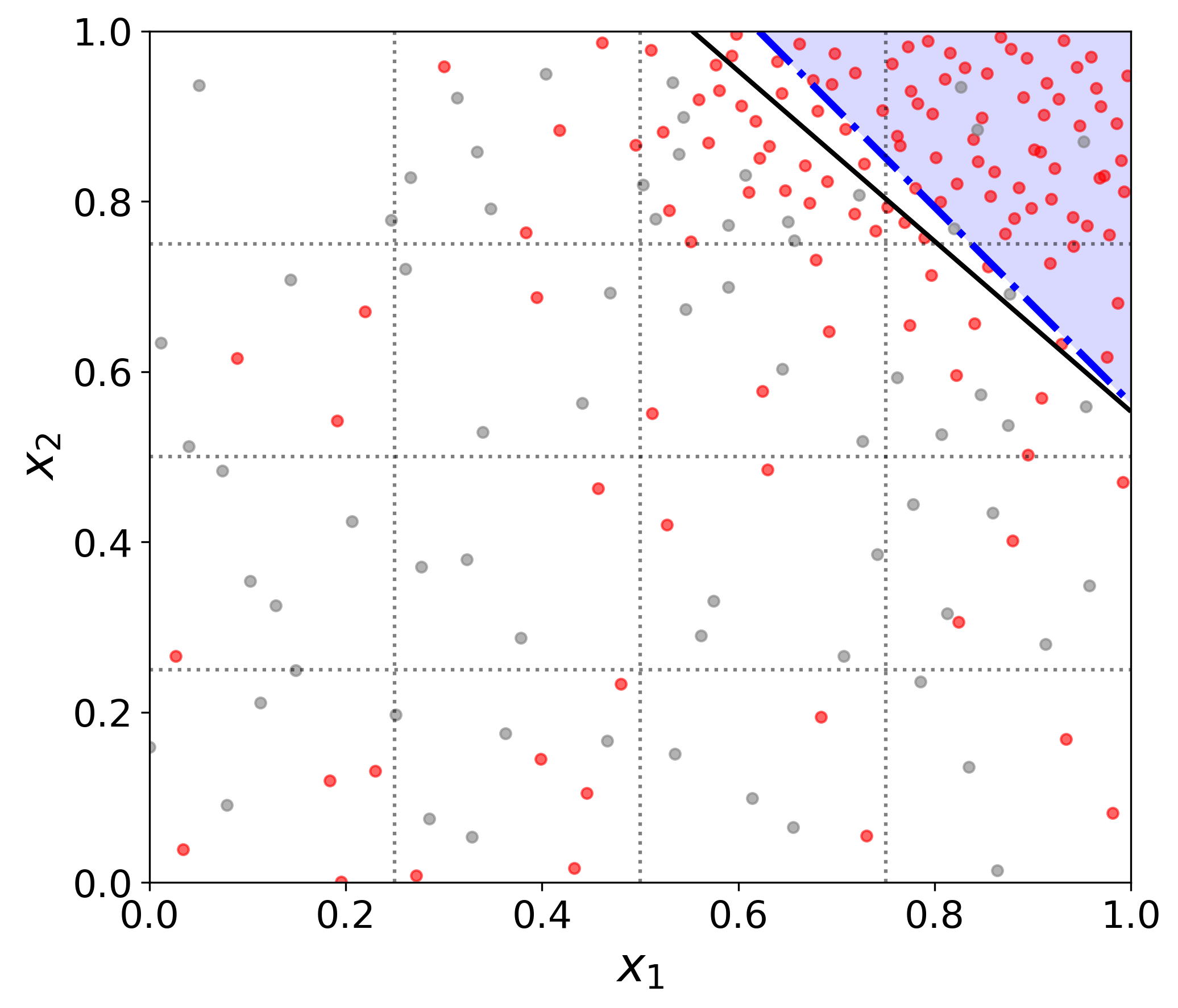}
    \end{subfigure}

 \caption{\it Visualization of the estimated reliability boundaries based on plug-in (first two columns with red dashed line) and the proposed method (last two columns with blue dot-dashed line) for $p_R=0.6$ (top panel) and $0.9$ (bottom panel) when the data is generated by {T-logit}. The training sets generated by SRS (first and third columns) and optimal designs (second and fourth columns) are displayed as scatter points on each plot. 
 The black solid line represents the true reliability boundary, and the light blue shaded area indicates the estimated reliable region yielded by different methods.}
\label{fig:est_comp}
\end{figure}

Figure~\ref{fig:est_comp} illustrates the estimated reliability boundaries drawn from the first trial among 500 independent replicates, illustrating both the estimation and design performance of the classical plug-in approach and our method. 
The proposed modeling-then-calibration approaches (the last two columns) deliver superior estimation performance. These calibrated methods successfully correct the severe boundary misalignments observed in the standard plug-in baselines (the first two columns). Specifically, they effectively mitigate both the dangerous overestimation of the reliability set at $p_R = 0.6$ and the severe underestimation at $p_R = 0.9$, ensuring that the predicted boundary tightly and safely aligns with the ground truth.
Another important observation from the third and fourth columns is that the proposed design places more experimental units in the true reliable region $\mathcal{X}_R$. This echoes the design principles we discussed in Section~\ref{sec:design}.
To quantify this more directly, we calculate the average number of design points that are allocated to  $\mathcal{X}_R$. 
We find that more than 110 
points are chosen in $\mathcal{X}_R$,
which means that more than 50\% of the total experimental budget ($n_{\rm train}$) can be used in the possible downstream analysis.

\subsection{Inference Performance under Risk Control}\label{subsec:CRC-sim}

Now we illustrate the inference performance under risk control, which is measured by the two kinds of risk, namely, \begin{align*}
    \label{eq:risk-decomp1}
     \text{Type~I~:}\quad  \risk_1(\lambda)& =     \E\bigl\{\bigl(p_R - Y\bigr) \, \mathbb{I}\bigl(\bm{x} \in \hat{\mathcal{X}}_R(\lambda), \bm{x} \notin {\mathcal{X}}_R\bigr) \bigr\}, \\
     \text{Type~II:}\quad   \risk_2(\lambda) &= \E\bigl\{ (Y-p_R) \mathbb{I}\bigl(\bm{x} \not\in \hat{\mathcal{X}}_R(\lambda), \bm{x} \in {\mathcal{X}}_R\bigr) 
        \bigr\} . 
    \end{align*}

\begin{table}[h]
    \centering\spacingset{1.0}
    \caption{\it Empirical median and interquartile range (IQR, in parentheses) Type I and Type II risk comparison across different design and analysis methods. 
    {Results that fail to control Type I risk at nominal level $\alpha(=0.01)$ are highlighted in red, and among those that succeed, the lowest Type II risk is highlighted in blue.}  }
    \label{tab:type_error_comp}
    \footnotesize 
    \setlength{\tabcolsep}{4.5pt} 
    \renewcommand{\arraystretch}{1.2} 

    \begin{tabular}{ c l l cccc cccc}
    \toprule
     \multirow{2}{*}{$\bm{p_R}$} & \multirow{2}{*}{\textbf{}} & \multirow{2}{*}{\textbf{}} & \multicolumn{4}{c}{\textbf{Plug-in}}& \multicolumn{4}{c}{\textbf{Ours}}\\
    \cmidrule(lr){4-7} \cmidrule(lr){8-11}
     & & & \multicolumn{2}{c}{SRS} & \multicolumn{2}{c}{Opt} & \multicolumn{2}{c}{SRS} & \multicolumn{2}{c}{Opt} \\
    \midrule
    
     \multirow{4}{*}{0.6} 
    & \multirow{2}{*}{T-logit} & Type I & 0.0018 & (0.0080) & {\red 0.0206} & (0.0173) & 0.0013 & (0.0063) & 0.0005 & (0.0041) \\ 
     & & Type II & 0.0012 & (0.0067) & 0.0000 & (0.0000) & 0.0020 & (0.0064) & {\blue 0.0009} & (0.0042) \\
     & \multirow{2}{*}{PWC} & Type I & {\red 0.0142} & (0.0134) & {\red 0.0366} & (0.0130) & 0.0029 & (0.0084) & 0.0024 & (0.0075) \\ 
     & & Type II & 0.0000 & (0.0000) & 0.0000 & (0.0000) & 0.0019 & (0.0060) & {\blue 0.0014} & (0.0051) \\ 
    \cmidrule(lr){1-11}
     \multirow{4}{*}{0.9} 
     &\multirow{2}{*}{T-logit} & Type I & 0.0000 & (0.0000) & 0.0000 & (0.0000) & 0.0006 & (0.0037) & 0.0004 & (0.0033) \\ 
     & & Type II & 0.0046 & (0.0010) & 0.0043 & (0.0016) & 0.0003 & (0.0016) & {\blue 0.0000} & (0.0011) \\
     & \multirow{2}{*}{PWC} & Type I & 0.0000 & (0.0000) & 0.0000 & (0.0000) & 0.0003 & (0.0037) & 0.0001 & (0.0032) \\ 
     & & Type II & 0.0063 & (0.0010) & 0.0055 & (0.0028) & 0.0013 & (0.0026) & {\blue 0.0011} & (0.0025) \\ 
    \bottomrule
    \end{tabular}
\end{table}

Table~\ref{tab:type_error_comp} compares the  Type I and Type II risks across different methods, where the nominal risk  $\alpha$ is chosen as  $0.01$. All the experiment setups are the same as Section~\ref{subsect:est-sim} except that $\lambda$ is chosen as described in Section~\ref{sec:conformal} (or \eqref{eq:hat-lambda1} more precisely). The two kinds of risk are evaluated empirically via $10000$ Monte Carlo points. 
As highlighted in the results, the plug-in estimators frequently violate the nominal risk level ($\alpha=0.01$), whereas our conformal risk control (CRC) successfully constrains the Type I risk below the target threshold. These empirical results support the theoretical guarantee established in
Theorem~\ref{thm:validity1}. Furthermore, in those cases where the plug-in
estimators successfully control the nominal risk, our method attains the smallest Type~II risk. This is consistent with Theorem~\ref{thm:lower_bound},
which shows that the proposed procedure is not overly conservative.

Another important application of the proposed risk control framework is in transfer learning settings, where researchers possess prior knowledge or historical models that are potentially misspecified but still contain useful information. Given a limited budget for new experiments, the goal is to construct a reliability set with guaranteed Type~I risk control, leveraging both the prior knowledge and the new observations. In such scenarios, researchers have the following three choices: 
\begin{itemize}
    \item \textbf{Naive}. Trust the prior models directly without any calibration, and output their implied reliability sets.
    \item \textbf{Plug-in}. Discard the prior entirely, conduct the new experiments, and estimate the reliability set solely from the new data using a standard plug‑in estimator.
    \item \textbf{Ours}. Treat the prior as approximately correct, use the same new experiments as a calibration set to select a threshold that ensures Type~I risk control, and thereby obtain a reliable reliability set with improved efficiency.
\end{itemize}
 
To illustrate, we consider a simulation study with target probability \(p_R = 0.9\) and the true data‑generating process follows a \(T\)-logit model. 
Two misspecified prior models, denoted \(\mathfrak M_1\) and \(\mathfrak M_2\), are assumed available. 
Both are linear logistic regression models with the same slope vector \((15, 15)^\T\) but different intercepts: \(-3.1\) for \(\mathfrak M_1\) and \(-27.1\) for \(\mathfrak M_2\). 
Only 100 new experimental trials are allowed. The Type~I risk is controlled at the nominal level \(\alpha = 0.01\).
The results are summarized in Table~\ref{tab:type_error_comp_robust}.

\begin{table}[h]
    \centering\spacingset{1.0}
    \caption{\it Empirical median and interquartile range (IQR, in parentheses) for Type I/II risk and volume of the symmetric difference (S-diff) across different analysis methods under $\mathfrak M_1$ and $\mathfrak M_2$. 
    {Results that fail to control Type I risk at the nominal level $\alpha=0.01$ are highlighted in red, and among those that succeed, the lowest Type II risk is highlighted in blue. The best performances of S-diff are highlighted in boldface.} }
    \label{tab:type_error_comp_robust}
    \footnotesize 
    \setlength{\tabcolsep}{8pt} 
    \renewcommand{\arraystretch}{1.2} 

    \begin{tabular}{c l cc cc cc}
    \toprule
     & \textbf{ } & \multicolumn{2}{c}{\textbf{Naive}} & \multicolumn{2}{c}{\textbf{Plug-in}} & \multicolumn{2}{c}{\textbf{Ours}} \\
    \midrule
    
    \multirow{3}{*}{$\mathfrak M_1$} 
    & Type I~  & {\red 0.4308} & (0.0057) & 0.0000 & (0.0000) & 0.0005 & (0.0009) \\ 
    & Type II & 0.0000 & (0.0000) & 0.0046 & (0.0010) & {\blue 0.0000} & (0.0000) \\
    & S-diff~~  & 0.8377 & (0.0011) & 0.0999 & (0.0005) & \textbf{0.0177} & (0.0104) \\ 
    \midrule
    
    \multirow{3}{*}{$\mathfrak M_2$} 
    & Type I~  & 0.0000 & (0.0000) & 0.0000 & (0.0000) & 0.0002 & (0.0005) \\ 
    & Type II & 0.0046 & (0.0010) & 0.0046 & (0.0010) & {\blue 0.0000} & (0.0000) \\ 
    & S-diff~~ & 0.0987 & (0.0005) & 0.0999 & (0.0005) & \textbf{0.0339} & (0.0053) \\
    \bottomrule
    \end{tabular}
\end{table}

Table~\ref{tab:type_error_comp_robust} reports the empirical Type I/II risks and the median of the volume of the symmetric difference (S‑diff) for the three methods under \(\mathfrak M_1\) and \(\mathfrak M_2\). The poor performance of the Naive method—severe Type I inflation under \(\mathfrak M_1\) and non‑negligible S‑diff under both priors—reflects the inherent misspecification of the prior models. 
Turning to the two data‑driven strategies, both the Plug‑in and our method successfully control the Type I risk at the nominal level. However, our proposed method achieves a substantially lower Type II risk (zero in both cases) and markedly smaller S‑diff compared to the plug‑in estimator. 
This demonstrates that by calibrating rather than discarding the misspecified prior, our framework not only guarantees risk control but also delivers higher estimation efficiency given the same experimental budget.

\vspace{-0.5cm}
\section{Design and inference for the car cash experiment}\label{sec:real-data}

To evaluate our framework in a practical setting, we analyze a vehicle safety experiment studied in \cite{imberg2025active}. 
The original data are collected from a simulator built by reconstructing the pre-crash kinematics of 44 rear-end collisions from a Volvo Car Corporation database in Sweden. 
The design space spans two variables: off-road glance duration and maximal deceleration. In the original study, glance duration is discretized into 67 equally spaced levels and
maximal deceleration into 15 equally spaced levels, yielding a full factorial grid of $67 \times 15 = 1{,}005$ design points at which the simulator is evaluated. 
For our analysis, both variables are normalized to the unit square $[0,1]^2$. Vehicle crash assessments inherently demand stringent safety standards; we therefore focus on the high-reliability threshold $p_R = 0.9$.

\begin{figure}[ht]
    \centering\spacingset{1.0}
    \begin{subfigure}[b]{0.3\textwidth}
        \centering
        \includegraphics[width=0.8\textwidth]{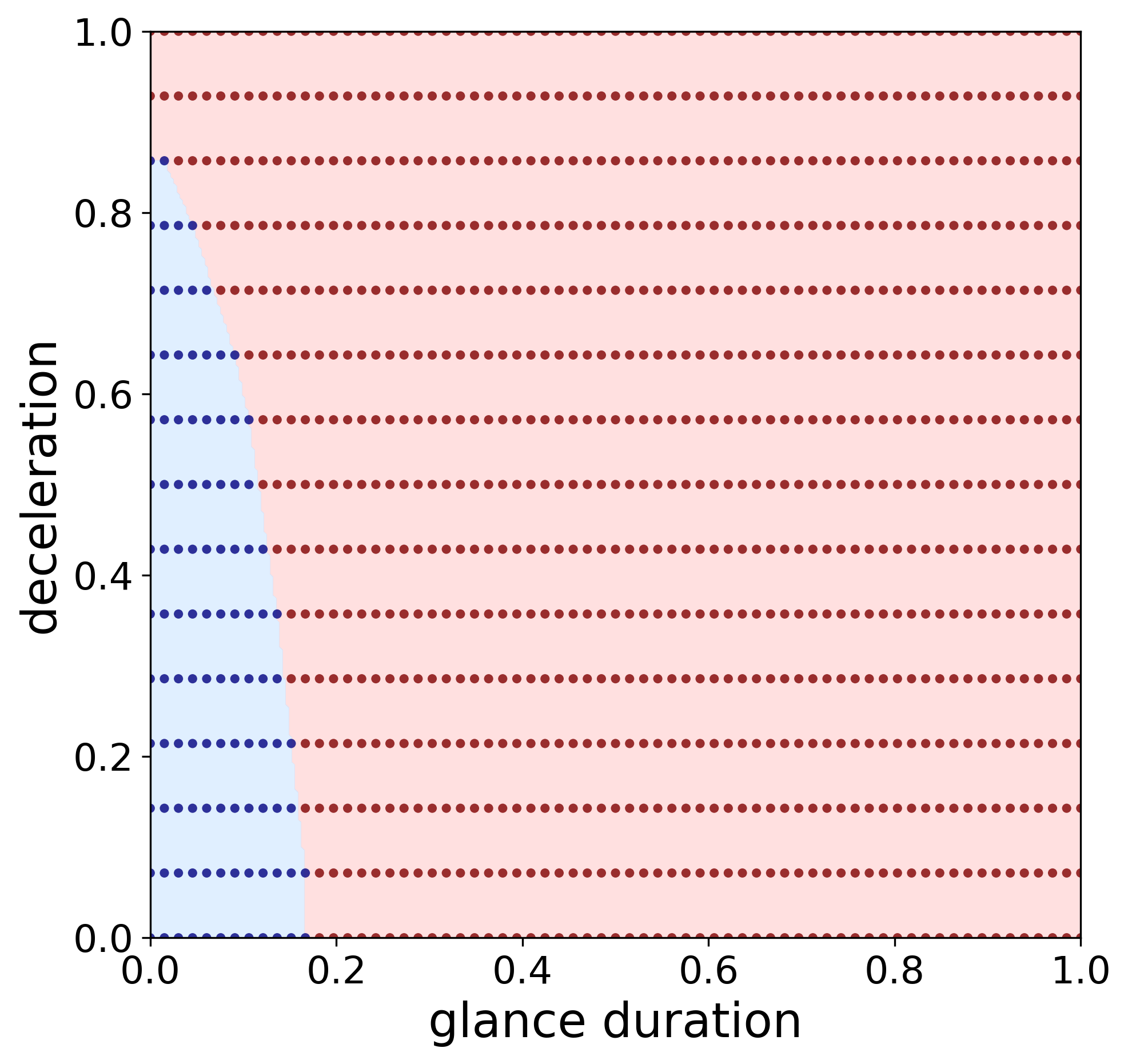}
    \end{subfigure}\hfill
    \begin{subfigure}[b]{0.3\textwidth}
        \centering
        \includegraphics[width=0.8\textwidth]{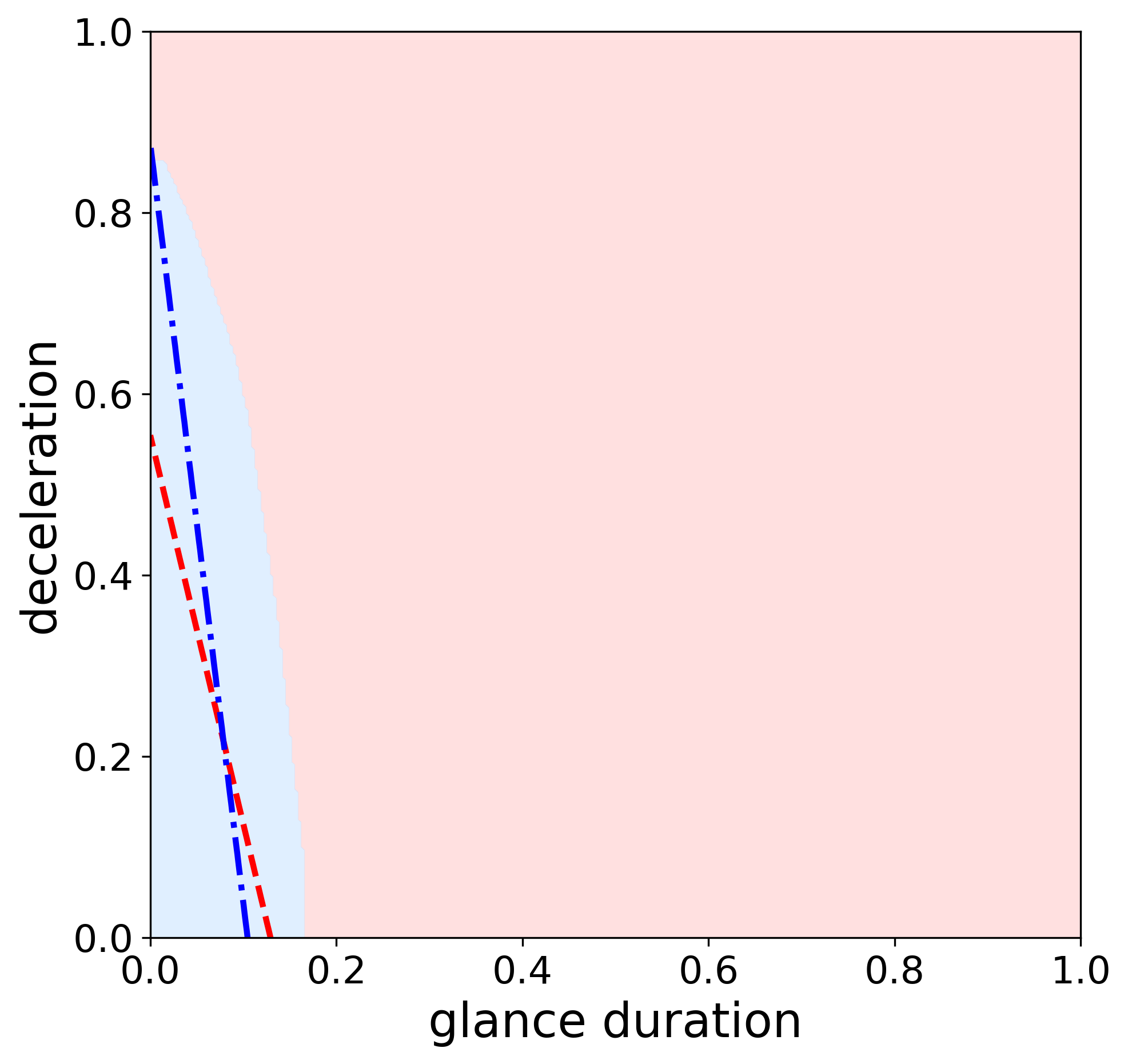}
    \end{subfigure}\hfill
    \begin{subfigure}[b]{0.3\textwidth}
        \centering
        \includegraphics[width=0.8\textwidth]{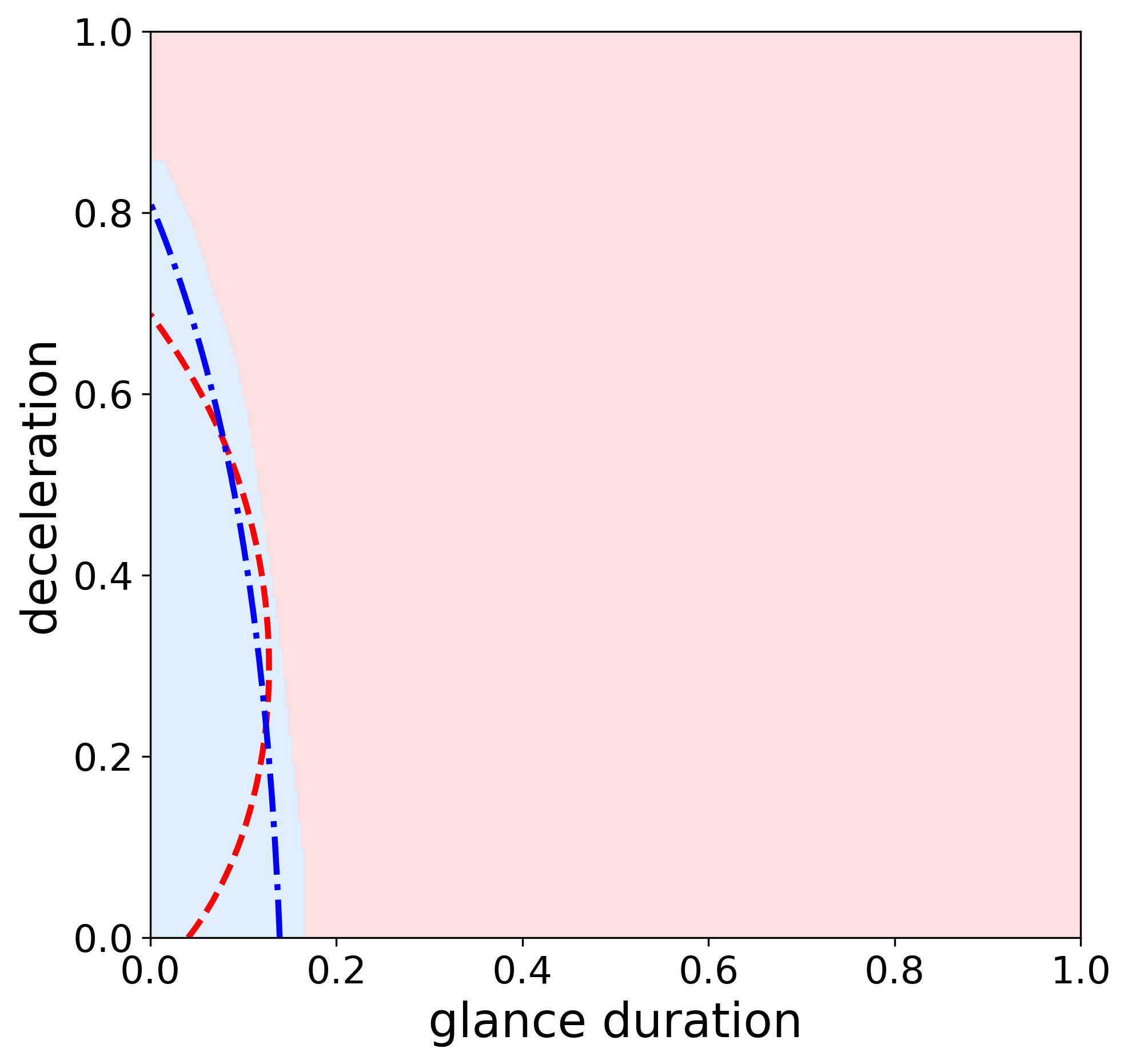}
    \end{subfigure}
    
    \caption{\it Left panel: 1005 raw data points, blue ($Y=1$) and red ($Y=0$). Middle and right panels: estimated boundaries from one simulation (using $n_{\train}=200$ training points and $n_{\cali}=100$ calibration points), with linear SVM and kernel SVM as working models, respectively. In both panels, the boundary obtained by the method proposed in this paper is represented by the blue dot-dashed line, while the red dashed line corresponds to the classical plug‑in approach.}
    \label{fig:real_data_boundaries}
\end{figure}

Because the simulator is deterministic, identical inputs always lead to identical outcomes. 
We construct a continuous target safety region as the convex hull of all safe instances ($Y=1$); any point outside this hull is strictly treated as a crash response ($Y=0$). 
The left panel of Figure~\ref{fig:real_data_boundaries}
visualizes the dataset with $Y=1$ marked in blue and $Y=0$ in red. 
The resulting decision boundary is sharp and irregular, making it challenging to capture adequately for smooth parametric models . Since the purely deterministic responses can lead to perfect separation and
 the maximum likelihood estimator of logistic regression fails, we opt to use linear and kernel support
vector machines (SVMs) as the working model $\pi$.

The experimental budget and parameter settings match those in Section~\ref{sec:numerical}: $n_{\train}=200$ training points and $n_{\cali}=100$ calibration points. 
 For each of the 500 independent repetitions, the design algorithm generates training points directly in the continuous design space $[0,1]^2$, and an independent calibration set is drawn uniformly from the same space. Because the simulator is
deterministic, the binary response for any point is uniquely determined by whether it lies inside the convex hull of the safe instances.

\begin{table}[htbp]
    \centering\spacingset{1.0}
    \caption{\it Empirical median and interquartile range (IQR, in parentheses) of the volume of the symmetric difference (S-diff) across different designs and analysis methods with different working models. As no optimal design strategies for SVMs within the plug-in estimator framework exist to our knowledge, the plug-in approach under optimal design is marked as \texttt{NA}.}
    \label{tab:est_comp_consistent}
    \footnotesize 
    \setlength{\tabcolsep}{3pt} 
    \renewcommand{\arraystretch}{1.2} 

    \begin{tabular}{ c cccc cccc}
    \toprule
       \multirow{2}{*}{\textbf{Classifier}} & \multicolumn{4}{c}{\textbf{Plug-in}}& \multicolumn{4}{c}{\textbf{Ours}}\\
    \cmidrule(lr){2-5} \cmidrule(lr){6-9}
    & \multicolumn{2}{c}{SRS} & \multicolumn{2}{c}{Opt} & \multicolumn{2}{c}{SRS} & \multicolumn{2}{c}{Opt} \\
    \midrule
    
   
     linear SVM & 0.0667 & (0.0106) & ~~~~\texttt{NA} & (\texttt{NA})~~~~ & 0.0280 & (0.0186) & \textbf{0.0226} & (0.0175) \\ 
     kernel SVM & 0.0354 & (0.0171) & ~~~~\texttt{NA} & (\texttt{NA})~~~~ & 0.0171 & (0.0124) & \textbf{0.0118} & (0.0104) \\ 
    \bottomrule
    \end{tabular}
\end{table}

These results align with our simulation findings: the proposed active design via the proposed analysis procedure substantially outperforms the classical plug‑in estimator with a random sampling design. 
For linear SVM, the S‑diff drops from 0.0667 to 0.0226 (a 66\% reduction); for kernel SVM, from 0.0354 to 0.0118 (a 67\% reduction). 
Moreover, kernel SVM consistently yields lower estimation error than linear SVM under both designs, reflecting the influence of the working model $\pi$'s quality. 
Even though kernel SVM is approximately correct, misspecification inevitably arises with finite samples, and thus our method still achieves a marked improvement over the baseline. The middle and right panel of Figure~\ref{fig:real_data_boundaries}  show the boundary obtained by our method (blue dot-dashed line) and the classical plug-in (red dashed line) with two different working models under $p_R=0.9$ and further illustrate the advantages of the new method. 

Parallel to Section~\ref{subsec:CRC-sim}, we also studied conformal risk control with the type-I risk being controlled at $\alpha=0.01$.
As expected, our proposed method successfully controls the type-I risk at the nominal level and also achieves the smallest type-II risk, which reflects that the Volvo car safety system performs better than we initially anticipated. Specifically, compared with the reliability set obtained by the plug‑in method, the true reliability set with $p_R = 0.9$ and Type~I risk controlled at $0.01$ tolerates a slightly higher deceleration and a slightly shorter glance duration.
Therefore, the vehicle can remain safe under more aggressive braking and shorter driver reaction times than the plug‑in estimator would suggest.

\begin{table}[h]
    \centering\spacingset{1.0}
    \caption{\it Empirical median and interquartile range (IQR, in parentheses) of Type I and Type II risk across different design and analysis methods with $\alpha=0.01$. 
    }
    \label{tab:error_comp_consistent}
    \footnotesize 
    \setlength{\tabcolsep}{3pt} 
    \renewcommand{\arraystretch}{1.2}

    \begin{tabular}{ l l cccccccc }
    \toprule
    \multirow{2}{*}{\textbf{Classifier}} & \multirow{2}{*}{} & \multicolumn{4}{c}{\textbf{Plug-in}} & \multicolumn{4}{c}{\textbf{Ours}} \\
    \cmidrule(lr){3-6} \cmidrule(lr){7-10}
    & & \multicolumn{2}{c}{SRS} & \multicolumn{2}{c}{Opt} & \multicolumn{2}{c}{SRS} & \multicolumn{2}{c}{Opt} \\
    \midrule
    
    \multirow{2}{*}{linear SVM} 
    & Type I   & 0.0000 & (0.0000) & ~~~\texttt{NA} & (\texttt{NA})~~~ & {0.0024} & (0.0055) & {0.0027} & (0.0053) \\ 
    & Type II & 0.0067 & (0.0011) & ~~~\texttt{NA} & (\texttt{NA})~~~ & 0.0024 & (0.0024) & {\blue 0.0019} & (0.0023) \\ 
    \midrule
    
    \multirow{2}{*}{kernel SVM} 
    & Type I   & 0.0000 & (0.0000) & ~~~\texttt{NA} & (\texttt{NA})~~~ & 0.0006 & (0.0032) & {0.0001} & (0.0015) \\ 
    & Type II & 0.0035 & (0.0017) & ~~~\texttt{NA} & (\texttt{NA})~~~ & 0.0014 & (0.0015) & {\blue 0.0010} & (0.0013) \\ 
    \bottomrule
    \end{tabular}
\end{table}


\section{Conclusion}

This paper developed a unified modeling-then-calibration framework for making reliability set estimation trustworthy under an imperfect working model. Rather than assuming the model to be correct, we accept that it is only useful, and we address the resulting trust deficit in two integrated stages.

The modeling stage is strengthened by an adaptive design (Section~\ref{sec:design}) that concentrates observations on a critical region covering the reliability set and its boundary, where estimation errors have the greatest impact on the final decision. This achieves space-filling over the region of interest while explicitly controlling the budget spent elsewhere, systematically improving model quality where it matters most (Theorem~\ref{thm:error_count_bound}).

The calibration stage employs a suite of procedures, each serving a distinct objective.
Risk minimizing calibration (Section~\ref{sec:method}) selects a threshold that corrects systematic decision errors; we show that it never increases risk, recovers the true set under mild separation conditions, and converges at a fast parametric rate (Theorems~\ref{lem:consistency}-\ref{thm:calibration}).
Conformal risk control (Section~\ref{sec:conformal}) provides a finite-sample guarantee on the risk of false inclusion, regardless of working model quality. 
Because the true set membership is unobservable, this guarantee requires a novel asymmetric surrogate loss, which provides control of the false inclusion risk (Theorems~\ref{thm:validity}-\ref{thm:validity1}).
Crucially, the results also hold when the working model is fixed in advance, as calibration alone can still deliver bias correction and finite-sample risk control (as illustrated in Section~\ref{subsec:CRC-sim}).

Beyond reliability set estimation, the unobservability of the true set is a recurring obstacle in set estimation problems with latent ground truth, where standard loss functions cannot be evaluated directly. Our surrogate loss construction offers a template for extending conformal risk control to such settings. A further design implication is that, when budgets are limited and model information is insufficient, a small calibration set is not an optional extra but a necessary component for risk control.

\bibliographystyle{agsm}
\spacingset{1.25}
\bibliography{reference}
\end{document}